\pdfoutput=1
\documentclass{amsart}

\usepackage{amsmath}
\usepackage{amsfonts}
\usepackage{cellspace}
\usepackage{longtable}
\newtheorem{thm}{Theorem}[section]
\newtheorem{lem}[thm]{Lemma}

\newtheorem{prop}[thm]{Proposition}

\newtheorem{defi}[thm]{Definition}

\theoremstyle{definition}
\newtheorem{rek}[thm]{Remark}
\newtheorem{exa}[thm]{Example}

\newtheorem{prob}{Problem}

\newcommand{\C}{\ensuremath{\mathbb{C}}}
\newcommand{\Z}{\ensuremath{\mathbb{Z}}}
\newcommand{\Q}{\mathbb{Q}}
\newcommand{\N}{\mathbb{N}}

\newcommand{\dx}{\mathrm{d}x}

\title[Factorial Ratios, Hypergeometric Series, and Step Functions]{Factorial Ratios, Hypergeometric Series, and a Family of Step Functions}
\author{Jonathan W. Bober}

\address{Department of Mathematics, University of Michigan, Ann Arbor, MI 48109, USA}
\email{bober@umich.edu}
\date{\today}

\begin{document}


\renewcommand{\theenumi}{(\roman{enumi})}
\renewcommand{\labelenumi}{\theenumi}

\renewcommand{\a}{\mathbf{a}}
\renewcommand{\b}{\mathbf{b}}
\renewcommand{\u}{\mathbf{u}}
\newcommand{\balpha}{\mathbf{\alpha}}
\newcommand{\bbeta}{\mathbf{\beta}}

\newcommand{\f}[5]{ { }_{#1}F_{#2}\left(\begin{array}{c}#3 \\ #4\end{array};#5\right)}
\newcommand{\lcm}{\mathrm{lcm}}
\newcommand{\floor}[1]{\left\lfloor#1\right\rfloor}
\newcommand{\intersect}{\cap}
\newcommand{\frc}[1]{\left\{#1\right\}}

\newcommand{\uab}{u_n(\a,\b)}
\newcommand{\uabfrac}{\frac{(a_1n)!(a_2n)!\cdots(a_Kn)!}{(b_1n)!(b_2n)!\cdots(b_Ln)!}}
\newcommand{\uabfracone}{\frac{(a_1n)!(a_2n)!\cdots(a_Kn)!}{(b_1n)!(b_2n)!\cdots(b_{K+1}n)!}}
\newcommand{\uabnfrac}[1]{\frac{(a_1(#1))!(a_2(#1))!\cdots(a_K(#1))!}{(b_1(#1))!(b_2(#1))!\cdots(b_L(#1))!}}
\newcommand{\fab}{f(x;\a,\b)}
\newcommand{\gab}{g(n;\a,\b)}
\newcommand{\fabx}[1]{f(#1;\a,\b)}
\newcommand{\fabsum}{\sum_{k=1}^K \floor {a_kx} - \sum_{l=1}^L \floor{b_lx}}
\newcommand{\fabaltsum}{\sum_{l=1}^L \frc{b_l x} - \sum_{k=1}^K \frc{a_k x}}
\newcommand{\sequence}[2]{#1_{1},#1_{2},\ldots #1_{#2}}

\newcommand{\efrac}[2]{e\left(\frac{#1}{#2}\right)}

\newcommand{\printlineno}{\hspace*{-.5in}\texttt{SOURCE LINE \#\the\inputlineno}}
\renewcommand{\printlineno}{}

\begin{abstract}
We give a complete classification of a certain family of step functions related to
the Nyman--Beurling approach to the Riemann hypothesis and previously studied by
V. I. Vasyunin. Equivalently, we completely describe when certain sequences of
ratios of factorial products are always integral.
Essentially,
once certain observations are made, this comes down to an application
of Beukers and Heckman's classification of the monodromy of the hypergeometric
function ${ }_nF_{n-1}$. We also note applications to the classification of cyclic
quotient singularities.
\end{abstract}

\maketitle
\printlineno

\section{Introduction}
In \cite{vasyunin-step-functions} V. I. Vasyunin considered the following problem, originating
from the Nyman--Beurling formulation of the Riemann hypothesis.
\begin{prob}\label{vasyunin-problem}\label{problem1}
Classify all step functions of the form
\begin{equation}\label{step-function-form}
    f(x) = \sum_{i=1} c_i \floor{\frac{x}{m_i}}, m_i \in \N, c_i \in \Z
\end{equation}
having the property that $f(x) \in \{0,1\}$ for all $x$.
\end{prob}

Vasyunin discovered some infinite families of solutions and fifty two additional
sporadic solutions and, based on the results of extensive computations,
formulated conjectures (\cite[Conjecture 8]{vasyunin-step-functions}, 
\cite[Conjecture 11]{vasyunin-step-functions}) that these lists were complete.

It follows easily from a theorem of Landau (see Proposition \ref{landau-lemma})
that this problem is equivalent to the following problem.


\begin{prob}\label{uab-problem}\label{problem2}
Let $\a \in \N^K, \b \in \N^{K+1}$, and set
\[
    \uab = \uabfracone.
\]
For what parameters $\a$ and $\b$ with $\sum a_i = \sum b_j$ is $\uab$ an integer for all $n$?
\end{prob}

It is immediately evident that in Problem \ref{vasyunin-problem} we may make a change
of variables $m_i' = Cm_i$ without changing the output of the function. Similarly, though
not as immediately obvious, $(\a,\b)$ is a solution to Problem \ref{uab-problem} if
and only if $(\frac{1}{d}\a, \frac{1}{d}\b)$ is a solution, where $d = \gcd(\a,\b)$. Thus, in
both cases it is enough to classify solutions with $\gcd 1$, as all other
solutions are multiples of these.

In \cite{R-V-factorials}, Rodriguez-Villegas observed that the generating function
attached to $\uab$, given by
\[
    u(\a,\b;z) = \sum_{n=0}^\infty \uab z^n,
\]
is in fact a hypergeometric series, and that $\uab$ is
always integral if and only if $u(\a,\b;z)$ is an algebraic function. This observation,
which is explained in Section \ref{hypergeometric-series-connection},
allows us to use the work of Beukers and Heckman \cite{BH-hypergeometric-function} to
give a complete answer to Problem \ref{uab-problem} and to prove Vasyunin's
conjectures about the completeness of his classification of such step functions.

Specifically, we prove the following theorem, which asserts that \cite[Conjecture 8]{vasyunin-step-functions} and
\cite[Conjecture 11]{vasyunin-step-functions} are true.

\printlineno

\begin{thm}\label{main-thm-cor}
Let
\[
    f(x) = \sum_{k=1}^N \floor{x \over m_k} - \sum_{k=N+1}^{2N + 1} \floor{x \over m_k},
\]
and suppose that $m_i \ne m_j$ for all $i \le N, j \ge N+1$,
and that
\[
\gcd(m_1, m_2, \ldots, m_{2N + 1}) = 1.
\]
Then $f(x)$ takes only the values $0$ and $1$ if and only if either
\begin{enumerate}
\item $f(x)$ takes one of the following forms:
\begin{equation}
    f(x) = \floor{x \over ab} - \floor{x \over b(a+b)} - \floor{x \over a(a+b)}
        \textrm{ where } \gcd(a,b) = 1,
\end{equation}
\begin{equation}
f(x) = \floor{x \over b(a-b)} + \floor{x \over 2a(a-b)} - \floor{x \over 2b(a-b)}
    - \floor{x \over a(a-b)} - \floor{x \over 2ab}
\end{equation}
where
\[
    \gcd(a,b) = \gcd(2,a-b) = 1 \textrm{ and } a > b > 0,
\]
\begin{equation}
f(x) = \floor{x \over {1 \over 2}b(a-b)} + \floor{x \over a(a-b)} - \floor{x \over b(a-b)}
    - \floor{x \over {1 \over 2}a(a-b)} - \floor{x \over ab}
\end{equation}
where
\[
    \gcd(a,b) = \gcd(2,a) = \gcd(2,b) = 1 \textrm{ and } a > b > 0,
\]
\begin{equation}
f(x) = \floor{x \over b(a+b)} + \floor{x \over a(a+b)} - \floor{x \over 2b(a+b)}
    -\floor{x \over 2a(a+b)} - \floor{x \over 2ab}
\end{equation}
where
\[
    \gcd(a,b) = \gcd(2,a+b) = 1,
\]
\begin{equation}
f(x) = \floor{x \over {1 \over 2}b(a+b)} + \floor{x \over {1 \over 2}a(a+b)} - \floor{x \over b(a+b)}
    -\floor{x \over a(a+b)} - \floor{x \over ab}
\end{equation}
where
\[
    \gcd(a,b) = \gcd(2,a) = \gcd(2,b) = 1,
\]
or
\item $f(x)$ is one of the $52$ sporadic step functions given by
\[
    (m_1, m_2, \ldots, m_N) = \left( \frac{M}{a_1}, \frac{M}{a_2}, \ldots, \frac{M}{a_N}\right)
\]
and
\[
    (m_{N+1}, m_{N+2}, \ldots, m_{2N + 1}) = \left( \frac{M}{b_1}, \frac{M}{b_2}, \ldots,
        \frac{M}{b_{N+1}}\right)
\]
for some $\a$ (or permutation of $\a$) and $\b$ (or permutation of $\b$)
listed in the second column of Table \ref{factorial-ratio-table},
where \[M = \lcm(a_1, a_2, \ldots, a_N, b_1, b_2, \ldots b_{N+1}).\]
\end{enumerate}
\end{thm}
This theorem is proved as a consequence of the equivalence of Problems \ref{problem1}
and \ref{problem2}, and the following theorem.

\begin{thm}\label{main-theorem}
Let
\[
    \uab = \uabfracone
\]
and suppose that $a_k \ne b_l$ for all $k,l$, that $\sum a_k = \sum b_l$,
and that
\[
\gcd(a_1, \ldots, a_K, b_1, \ldots, b_{K+1}) = 1.
\]
Then $\uab$ is an integer for all $n$ if and only if either
\begin{enumerate}
\item $u_n = \uab$ takes one of the following forms:
\begin{equation}
u_n = \frac{[(a+b)n]!}{(an)!(bn)!} \textrm{ for $\gcd(a,b) = 1$},
\end{equation}
\begin{equation}
u_n = \frac{(2an)!(bn)!}{(an)!(2bn)![(a-b)n]!} \textrm{ for $\gcd(a,b) = 1$ and $a > b$},
\end{equation}
\begin{equation}
u_n = \frac{(2an)!(2bn)!}{(an)!(bn)![(a+b)n]!} \textrm{ for $\gcd(a,b) = 1$}
\end{equation}
or
\item $(\a,\b)$ is one of the $52$ sporadic parameter sets listed
in the second column of Table $\ref{factorial-ratio-table}$.
\end{enumerate}
\end{thm}
In Section \ref{quotient-singularities} we note briefly that this Theorem has immediate applications
to the classification of cyclic quotient singularities.

It is a standard elementary exercise (see, for example, \cite[Section 4.5]{apostol-introduction-to-analytic-number-theory})
to use the integrality of factorial ratios such as these
to prove Chebyshev's elementary estimate for the prime counting function $\pi(x)$ --- namely,
that there exist numbers $c_1$ and $c_2$ with $c_1 < 1 < c_2$ such that
\[
    c_1 \frac{x}{\log x} \le \pi(x) \le c_2 \frac{x}{\log x}
\]
for all $x$ large enough.
It is interesting to note that the best such constants achievable from the factorial
ratio sequences listed in Theorem \ref{main-theorem} are those which were 
discovered by Chebyshev \cite{chebyshev-elementary-pnt-estimates}. The factorial ratio sequence
that gives these constants ($c_1 \approx .92$ and $c_2 \approx 1.11$) is
\begin{equation}\label{chebyshev-ratio}
    u_n = \frac{(30n)!n!}{(15n)!(10n)!(6n)!},
\end{equation}
which was in fact used by Chebyshev. (In truth, although Chebyshev's method is in some
sense equivalent to using factorial ratios, it is not quite the same. Chebyshev actually
makes more direct use of the fact that the corresponding step function
\[
    f(x) = \floor{x} - \floor{x \over 2} - \floor{x \over 3} - \floor{x \over 5} + \floor{x \over 30}
\]
takes on only the values $0$ and $1$.)
The monodromy group associated with this $u_n$ --- that is, the monodromy group of the hypergeometric 
differential equation satisfied by $\u(\a,\b)$ --- is $W(E_8)$, the Weyl group of the $E_8$ root system,
which is the largest ``sporadic'' finite primitive complex reflection group.

\subsection{Relation to Nyman--Beurling Formulation of the Riemann hypothesis}

Let 
\[
    \rho_\alpha(x) = \floor{\alpha \over x} - \alpha\floor{1 \over x} \in L^2( (0,1) )
\]
There are various similar statements of the Nyman--Beurling formulation
of the Riemann hypothesis. One such statement is

\begin{thm}[Nyman \cite{nyman-thesis}, B\'aez-Duarte \cite{MR2057270}]
The Riemann hypothesis is true if and only if
\[
    \mathrm{span}_{L^2(0,1)}\left\{\rho_\alpha(x) | {1 \over \alpha} \in \N\right\} = L^2( (0,1) ).
\]
\end{thm}
 
\begin{rek}
B. Nyman \cite{nyman-thesis} originally proved a weaker form of this theorem, with
$\alpha \in (0,1)$ instead of $1/\alpha \in \N$. The strengthening was conjectured for some
time and was proved by B\'aez-Duarte \cite{MR2057270}, however, B\'aez-Duarte works
with a slightly different formulation.
For a proof of the weaker version of this theorem in the above form, using theory
of Hardy spaces, and more discussion, a good
expository article is Balazard and Saias \cite{MR1759315}. An alternate, purely complex-analytic proof,
is given by Burnol \cite{MR2136671}.
\end{rek}

In fact, Beurling \cite{MR0070655} noticed that it is enough to approximate the constant function by linear combinations of
the functions $\rho_{1/n}(x)$. Specifically, one has

\begin{thm}[Beurling \cite{MR0070655}, B\'aez-Duarte \cite{MR2057270}]
The Riemann hypothesis is true if and only if
\[
    \chi_{(0,1)}(x) \in \mathrm{span}_{L^2(0,1)} \frc{\rho_\alpha(x) | {1 \over \alpha} \in \N},
\]
where $\chi_{(0,1)}(x)$ is the characteristic function of the interval $(0,1)$.
\end{thm}

This is Vasyunin's motivation for studying Problem \ref{vasyunin-problem}. Making the change of variables
$x \to 1/x$, finite linear combinations of functions $\rho_\alpha$ become functions of the form
\[
    \sum_{k=1}^n c_k \floor{x \over m_k}, m_k \in \N
\]
such that
\[
    \sum_{k=1}^n {c_k \over m_k} = 0,
\]
and now the function space is $L^2( (1,\infty), {\dx \over x^2})$. Thus, the Nyman--Beurling criterion
can be stated as: The Riemann hypothesis is equivalent to the existence of a sequence of functions
$\{\phi_n(x)\}$ all of the above form such that
$\phi_n(x) \rightarrow 1$ in $L^2( (1,\infty), {\dx \over x^2})$.
Vasyunin studies families of functions that take only the values $0$ and $1$ to construct
pointwise approximations to the constant function, and then explores the question stated in
Problem \ref{vasyunin-problem}. In particular, Vasyunin lists some infinite families of solutions
to Problem \ref{vasyunin-problem} and $52$ ``sporadic'' solutions, found by extensive computer search,
and he conjectures that his list is complete. Theorem \ref{main-thm-cor} says that this is indeed the case.
\footnote{The step function corresponding to line 12 of Table \ref{factorial-ratio-table} is not listed in
\cite{vasyunin-step-functions}. However, Vasyunin states that he found $21$ seven term step
functions, but lists only $20$, so its omission must be a transcription error.}

B\'aez-Duarte \cite{baez-duarte-vasyunin-correction} has shown, as Vasyunin
was well aware, that at least one of the pointwise approximations to
the constant function that Vasyunin constructs diverges in $L^2( (1,\infty), {\dx \over x^2})$.
It seems likely that the other sequences constructed by Vasyunin diverge as well. In light of
this, it seems that if this explicit approach is to give any insight into the
Riemann hypothesis, it will be necessary to study step functions that take on more than two values.
One natural question to ask might be: Can we classify all step functions of the
form \eqref{step-function-form} that take only values in the set $\{0,1,\dots,D\}$
for fixed $D$? For $D > 1$, this question seems difficult. Using substantially different
methods, a partial result on this question is obtained in a separate paper in preparation by J. Bell
and the author.

\subsection{Application: Cyclic quotient singularities}\label{quotient-singularities}

A. Borisov \cite{borisov-quotient-singularities} has noted that there is a connection
between integral factorial ratios and cyclic quotient singularities. In particular, we
note that Borisov shows that Vasyunin's conjecture (Theorem \ref{main-thm-cor}) implies the following.
(See \cite[Conjecture 1]{borisov-quotient-singularities}.)

\begin{thm}
Suppose $d \ge 5$ and we have a one-parameter family of Gorenstein cyclic quotient
singularities of dimension $2d + 1$ with Shukarov minimal log-discrepancy $d$. Then up
to the permutation of the coordinates in the $T^{(2d + 1)}$, the corresponding points lie
in the subtorus $x_1 + x_2 = 1$.
\end{thm}

For more on this connection, the reader should consult Borisov's paper.

\subsection{Notation:}
$(\alpha)_n$ denotes the rising factorial,
\[
    (\alpha)_n := (\alpha)(\alpha+1)(\alpha+2)\cdots(\alpha+n-1),
\]
and, for $\alpha = (\alpha_1, \ldots, \alpha_n)$, $\beta = (\beta_1, \ldots, \beta_m)$,
${ }_nF_m(\alpha; \beta; z)$ is the hypergeometric function
\[
    { }_nF_m(\alpha; \beta; z) = \sum_{k=0}^\infty
        { (\alpha_1)_k(\alpha_2)_k\cdots(\alpha_n)_k \over
            (\beta_1)_k(\beta_2)_k\cdots(\beta_m)_k}{z^k \over k!}.
\]
Also, $e(x) := \exp(2\pi i x) := e^{2\pi i x}$ and $\zeta_n = e(1/n)$ denotes the primitive
$n$th root of unity with smallest positive argument.

\subsection{Acknowledgements}
I am grateful to J.C. Lagarias for suggesting the problem considered in this
paper, and for encouragement and helpful comments
on earlier drafts of this paper. I was aided
in the discovery of some of these results by computer code I wrote for the
freely available SAGE Mathematics Software \cite{sage-2.6}, even though none
of these computations are actually needed to prove the given results.

\section{Some Preliminary Definitions and Notation}\label{lemmas-and-notations-section}
Throughout this paper, unless otherwise specified, $\a$ and $\b$ denote ordered
tuples of positive integers
\[
 \a = (a_1, a_2, \dots, a_K)
\]
and
\[
 \b = (b_1, b_2, \dots b_L).
\]
$\uab$ denotes the factorial ratio
\[
    \uab = \uabfrac.
\]
and $\fab$ denotes the step function
\[
    \fab = \sum_{k=1}^K \floor{a_kx} - \sum_{l=1}^L \floor{b_lx}
\]
Only later will we specify the requirement that $L = K+1$.

It is useful to attach certain polynomials to $\uab$ as follows.

\printlineno

\begin{defi}\label{P-Q-definition}
Define
$P(x) = P(\a,\b;x) \in \Z[x]$ and $Q(x) = Q(\a,\b;x) \in \Z[x]$ to be relatively prime polynomials
such that
\[
    \frac{P(x)}{Q(x)} = \frac{(x^{a_1} -1)(x^{a_2} - 1)\cdots(x^{a_K} - 1)}
        {(x^{b_1} - 1)(x^{b_2} - 1)\cdots(x^{b_L} - 1)}
\]
Then for some $\alpha_1 \le \alpha_2 \le \ldots \le \alpha_d$
and $\beta_1 \le \beta_2 \le \ldots \beta_d$, with $0 < \alpha_i,\beta_j \le 1$, $P$ and $Q$ factor
in $\C[x]$ as
\[
    P(x) = (x - e(\alpha_1))\cdots(x - e(\alpha_d))
\]
and
\[
    Q(x) = (x - e(\beta_1))\cdots(x - e(\beta_d)).
\]
where $e(x) = e(2\pi i x)$.

Set $\alpha(\a,\b) = \{\alpha_1, \ldots, \alpha_d\}$ and
$\beta(\a,\b) = \{\beta_1, \ldots, \beta_d\}$.
\end{defi}
\begin{rek}\label{misc-rek-one}
It is not too hard to see that for $x \in [0,1]$
\[
    \fab = \#\{\alpha_i | \alpha_i \le x\} - \#\{\beta_i | \beta_i \le x\}.
\]
Thus, an alternate definition for $\fab$ could be
\[
    \fab = \#(\alpha(\a,\b) \cap [0, \frc{x}]) - \#(\beta(\a,\b) \cap [0,\frc{x}])
\]
\end{rek}

We will occasionally make use of the notion of the interlacing of two sets,
so we state the following formally as a definition.

\begin{defi}[Interlacing]
 We say that two finite sets of real numbers $A$ and $B$ interlace
if the function
\[
 f(x) = \#\left( (-\infty, x) \cap A\right) - \#\left( (-\infty, x) \cap B\right)
\]
either takes only the values $0$ and $1$, or takes only the values $-1$ and $0$. In other
words, there is an element of $A$ in between any two elements of $B$, and an element of
$B$ in between any two elements of $A$.

We say that two sets $A$ and $B$ of complex numbers on the unit circle
\emph{interlace on the unit circle} if their arguments interlace on the
real line, where we take the argument of a complex number to be in $[0, 2\pi)$.
\end{defi}

\printlineno

\section{The connection between step functions and factorial ratios: The equivalence
of Problems \ref{vasyunin-problem} and \ref{uab-problem}}

In this section we will prove Theorem \ref{main-thm-cor} using Theorem \ref{main-theorem} and
the equivalence of Problems \ref{vasyunin-problem} and \ref{uab-problem}.
We begin by stating a rather general theorem of Landau \cite{landau-factorials} that connects the
integrality of factorial ratios with the nonnegativity of related step functions.

\begin{thm}[Landau \cite{landau-factorials}]\label{landau-theorem} Let $a_{k,s}, b_{l,s} \in \Z_{\ge 0}, 1 \le k \le K, 1 \le l \le L,
1 \le s \le r$ and let
\[
    A_k(x_1, x_2, \ldots, x_r) = \sum_{s=1}^r a_{k,s}x_s
\]
and
\[
    B_l(x_1, x_2, \ldots, x_r) = \sum_{s=1}^r b_{l,s}x_s.
\]
(That is, $A_k$ and $B_l$ are linear forms in $r$ variables with nonnegative integral
coefficients.) Then the factorial ratio
\[
    \frac{\prod_{k=1}^K A_k(x_1, x_2, \ldots, x_r)!}{\prod_{l=1}^L B_l(x_1, x_2, \ldots, x_r)!}
\]
is an integer for all $(x_1, \ldots, x_r) \in \Z_{\ge 0}^r$ if and only if the step function
\[
F(y_1, \ldots, y_r) = \sum_{k=1}^K \floor{A_k(y_1, \ldots, y_r)} - 
        \sum_{l=1}^L \floor{B_l(y_1, \ldots, y_r)}
\]
is nonnegative for all $(y_1, \ldots, y_r) \in [0,1]^r$.
\end{thm}
\begin{proof}
 See \cite{landau-factorials}.
\end{proof}

The special case of this that we will use is the following.

\begin{lem}\label{landau-lemma}
Let
\[
    u_n = \uab = \frac{(a_1n)!(a_2n)!\cdots(a_Kn)!}
               {(b_1n)!(b_2n)!\cdots(b_Ln)!}.
\]
Then $u_n$ is an integer for all $n$ if and
only if the function
\[
    f(x) = \fab = \sum_{k=1}^K \floor{a_kx} - \sum_{l=1}^L \floor{b_lx}
\]
is positive for all $x$ between $0$ and $1$. 
\end{lem}
\begin{proof}
 Take $A_k(x) = a_kx$ and $B_l(x) = b_lx$ in Theorem \ref{landau-theorem}.
\end{proof}

It also turns out that if $\fab$ is ever negative, then every prime that is large enough
occurs as a factor in the denominator of $\uab$ for some $n$.

\begin{lem}\label{uab-denominator-lemma}
Let
\[
    u_n = \uab = \frac{(a_1n)!(a_2n)!\cdots(a_Kn)!}
               {(b_1n)!(b_2n)!\cdots(b_Ln)!}.
\]
If $u_n$ is not
an integer for some $n$, then there exists some integer $P$ such that
for each prime $p > P$ there exists some $n$ such that $v_p(u_n) < 0$
(where $v_p(u_n)$ is the $p$-adic valuation of $u_n$).
\end{lem}

\begin{proof}
Consider the $p$-adic valuation of $n!$. We have
\[
    v_p(n!) = \sum_{\alpha = 1}^\infty \floor{\frac{n}{p^\alpha}}.
\]
Thus we have
\[
    v_p(u_n) = \sum_{\alpha = 1}^\infty f\left(\frac{n}{p^\alpha}\right),
\]
where $f(x) = \fab$

\printlineno

Assuming that $u_n$ is not always an integer, we know
from Lemma \ref{landau-lemma} that $f(x)$ is negative for some $x$. Since
$f$ is a step function, it follows that there is some interval, say
$[\beta, \beta + \epsilon]$ such that $f(x) < 0$ for all $x \in [\beta, \beta + \epsilon]$.
Additionally, we know that there is some $\delta > 0$ such that $f(x) = 0$ for
all $x \in [0,\delta]$. If we could find some $n$ and $p$ such that $n/p \in [\beta, \beta + \epsilon]$
and $n/p^2 \in [0,\delta]$, then we would have $f(n/p) < 0$ and $f(n/p^\alpha) = 0$ for
all $\alpha > 1$, and so we would clearly have $\nu_p(u_n) < 0$.

Now, such an $n$ and $p$ need to simultaneously satisfy the two inequalities
\[
    p\beta \le n \le p(\beta + \epsilon)
\]
and
\[
    0 \le n \le p^2\delta.
\]
For $p$ large enough, say $p > P_1$, we have $p^2\delta > p(\beta + \epsilon)$,
so it is sufficient for $n$ and $p$ to satisfy the first of the inequalities.
Moreover, for any $p$ large enough, say $p > P_2$, we have $p\epsilon > 1$, so that
there will in fact be an integer $n$ between $p\beta$ and $p(\beta + \epsilon)$. So
in fact, for any $p > P = \max(P_1, P_2)$ we have that there exists an $n$ such that
$\nu_p(u_n) < 0$.
\end{proof}

\printlineno

Along with Lemma \ref{landau-lemma}, the following lemma, which is a simple
generalization of \cite[Proposition 3]{vasyunin-step-functions} will yield the full equivalence of Problems
\ref{vasyunin-problem} and \ref{uab-problem}.

\begin{lem}\label{k-minus-l-condition}
Suppose that $f(x)$ is a function of the form
\[
    f(x) = \sum_{k=1}^K \floor{a_k x} - \sum_{l=1}^L \floor{b_l x}
\]
with $a_k, b_l \in \Z$ and that $f(x)$ is bounded. Then
$\sum_{k=1}^K a_k = \sum_{l=1}^L b_l$ and, for any $n$, there
exists some $x$ such that $f(x) = -n$ if and only if there exists
some $x'$ such that $f(x') = L-K + n$. In particular, $f(x)$ is nonnegative
if and only if the maximum value of $f$ is $L-K$.
\end{lem}
\begin{proof}
The first assertion is clear, for if $\sum a_k \ne \sum b_l$, then
$f(x)$ is unbounded. Now we know that $f(x)$ is periodic with period $1$. Now, for
any $z$ that is not an integer we have $\floor{z} + \floor{-z} = -1$, so for
any $z$ for which none of $a_i z$, $b_j z$ is an integer, we have
\[
    f(z) + f(-z) = L-K,
\]
from which the assertion follows.
\end{proof}

The following lemma describes explicitly the equivalence between Problems \ref{vasyunin-problem}
and \ref{uab-problem}, in a slightly more general form.

\begin{lem}\label{equivalence-lemma}
Let $\a = (a_1, a_2, \ldots, a_K), \b = (b_1, b_2, \ldots, b_L)$,
and put
\[M = \lcm(a_1, a_2, \ldots, a_K, b_1, b_2, \ldots, b_L).\]
Set
\[
(m_1, m_2, \ldots, m_{L+K}) = \left(\frac{M}{a_1}, \frac{M}{a_2}, \ldots, \frac{M}{a_L}, 
    \frac{M}{b_1}, \frac{M}{b_2}, \ldots, \frac{M}{b_k}\right).
\]
Then the following are equivalent:
\begin{enumerate}
\item
\[
    f(x) = \sum_{i=1}^K \floor{x \over m_i} - \sum_{i=K+1}^{K+L} \floor{x \over m_i}
\]
takes on values only in the range $0, 1, \ldots, L-K$.
\item $\sum_{k=1}^K a_k = \sum_{l=1}^L b_l$ and
\[
    u_n = \uabfrac
\]
is an integer for all $n \in \N$.
\end{enumerate}
Moreover, if $L-K = 1$, we may add
\begin{enumerate}
 \item[(iii)] $\alpha(\a,\b)$ and $\beta(\a,\b)$ interlace.
\end{enumerate}

\end{lem}
\begin{proof}
$f(x)$ differs from $f(x;\a,\b)$ only by a change of variables, so
Proposition \ref{landau-lemma} tells us that $u_n \in \Z$ for all $n$ if and only
if $f(x) > 0$ for all $x \in [0,1]$. Additionally, the boundedness of $f(x)$ is
equivalent to the statement that $\sum a_k = \sum b_l$, and Lemma \ref{k-minus-l-condition}
tells us that $f(x)$ is bounded and nonnegative if and only if it maximum value is
$L-K$.

Part (iii) is simply an issue of terminology, and follows directly from the alternate definition
of $\fab$ given in Remark \ref{misc-rek-one}.
\end{proof}

Using this equivalence, we can prove Theorem \ref{main-thm-cor} as an application
of Theorem \ref{main-theorem}.

\begin{proof}[Proof of Theorem~{\rm\ref{main-thm-cor}} (using Theorem~{\rm\ref{main-theorem}})]
The only complication that remains is that of classifying solutions with greatest
common divisor $1$. Consider the map $\phi : \N^K \times N^L \rightarrow \N^K \times \N^L$ given by
\[
\phi(a_1, a_2, \ldots, a_K, b_1, \ldots, b_L) = 
  \left({M \over a_1}, {M \over a_2}, \ldots, {M \over a_K}, {M \over b_1}, \ldots, {M \over b_L}\right)
\]
where
\[
    M = \lcm(a_1, a_2, \ldots, a_K, b_1, \ldots, b_L).
\]
The image of $\phi$ is all $(K + L)$-tuples with greatest common divisor $1$ and $\phi$ is
bijective on this subset. Thus $\phi$, in combination with Lemma \ref{equivalence-lemma},
gives a bijection between solutions to Problem
\ref{vasyunin-problem} with greatest common divisor $1$ and solutions to Problem \ref{uab-problem}
with greatest common divisor $1$.

When we apply this map to the three families of solutions to Problem \ref{uab-problem}
listed in Theorem \ref{main-theorem}, we get the five families of solutions
to Problem \ref{vasyunin-problem} listed in Theorem \ref{main-thm-cor}, and the $52$
sporadic solutions are given by the $52$ sporadic solutions to Problem \ref{uab-problem}.
\end{proof}

\section{The connection between factorial ratios and hypergeometric series}\label{hypergeometric-series-connection}

Rodriguez-Villagas \cite{R-V-factorials} observed a connection between hypergeometric
series and factorial ratio sequences. The purpose of this section is to write out this connection
explicitly in order to use it for our classification.

We begin with a lemma to show that the generating function for $\uab$ is in
fact a hypergeometric series.
\begin{lem}
Let
\[
    \uab = \uabfrac
\]
and
\[
    \u(\a,\b;z) = \sum_{n=0}^\infty \uab z^n.
\]
Let $\alpha = \alpha(\a,\b)$ and $\beta = \beta(\a,\b)$. If $L > K$, set
$d' = d, \alpha' = \alpha$, and $\beta' = \{\beta_1, \ldots, \beta_{d-1}\}$.
Otherwise, set $d' = d+1, \alpha' = \{\alpha_1, \ldots, \alpha_d, 1\}$, and
$\beta' = \beta$.

Then $\u(\a,\b)$ is the hypergeometric series
\[
    \u(\a,\b;z) = { }_{d'}F_{d'-1}\left(\alpha';\beta';Cz\right)
\]
where
\[
    C = \frac{a_1^{a_1}\cdots a_K^{a_K}}{b_1^{b_1}\cdots b_L^{b_L}}.
\]
\end{lem}

\printlineno

\begin{proof}
Examine the ratio between two consecutive terms
\begin{eqnarray*}
    A(n+1)  &=& \frac{u_{n+1}(\a,\b)}{\uab} \\ &=&\uabnfrac{n+1} \\ && \times \left[\uabfrac\right]^{-1}
\end{eqnarray*}
After cancellation, this can be written as
\[
    A(n+1) = \frac{(a_1n + 1) (a_1n + 2) \cdots (a_1n + a_1)(a_2n + 1)\cdots(a_Kn + a_K)}
        {(b_1n + 1)(b_1n + 2)\cdots(b_1n + b_1)(b_2n + 1)\cdots(b_Ln + b_L)}.
\]
Now if we factor out the coefficients of $n$ in each term we get
\[
    A(n+1) = C
        \frac{\left(n + \frac{1}{a_1}\right)\left(n + \frac{2}{a_1}\right)\cdots\left(n + \frac{a_1}{a_1}\right)\left(n + \frac{1}{a_2}\right)\cdots\left(n + \frac{a_K}{a_K}\right)}
             {\left(n + {1 \over b_1}\right)\left(n + {2 \over b_1}\right)\cdots\left(n + {b_1 \over b_1}\right)
                \left(n + {1 \over b_2}\right)\cdots\left(n + {b_L \over b_L}\right)},
\]
where
\[
    C = \frac{a_1^{a_1}\cdots a_K^{a_K}}{b_1^{b_1}\cdots b_L^{b_L}}.
\]
If we remove the common factors in the fraction then for exactly the same $\alpha$ and
$\beta$ as in Definition \ref{P-Q-definition} we have
\[
    A(n+1) = C\frac{(n + \alpha_1)\cdots(n + \alpha_d)}{(n + \beta_1)\cdots(n + \beta_d)}.
\]
Now, $u_0(\a,\b) = 1$, so we have in general
\[
    \uab = \prod_{k=1}^n A(k) = C^n\frac{(\alpha_1)_n\cdots(\alpha_{d})_n}{(\beta_1)_n\cdots(\beta_{d})_n},
\]
Now, if $\beta_{d} = 1$, we set $d' = d$. Otherwise, set $d' = d + 1$ and $\alpha_{d+1} = 1$. Either
way we have
\[
    \uab = \frac{C^n}{n!}\frac{(\alpha_1)_n\cdots(\alpha_{d'})_n}{(\beta_1)_n\cdots(\beta_{d'-1})_n}
\]
and so
\[
    \u(\a,\b;z) = \sum_{n=0}^\infty\frac{(\alpha_1)_n\cdots(\alpha_{d'})_n}{(\beta_1)_n\cdots(\beta_{d'-1})_n}
        \frac{(Cz)^n}{n!} = { }_{d}F_{d' - 1}(\alpha,\beta;Cz).
\]
To finish, we notice that $\beta_d = 1$ if and only if $L > K$.
\end{proof}

\begin{exa}
Let $\a = (30,1)$ and $\b = (15,10,6)$ and
\[
    u_n = \uab = \frac{(30n)!n!}{(15n)!(10n)!(6n)!}.
\]
Consider the ratio $\frac{u_{n+1}}{u_n}$. This is
\[
\frac{ (30n + 1)(30n + 2)\ldots(30n + 30)(n+1)}{(15n + 1)\ldots(15n + 15)(10n + 1)\ldots(10n + 10)(6n + 1)\ldots(6n+6)}.
\]
Factoring out the coefficients of $n$ in each term in the products, we get
\[
\frac{30^{30} (n + \frac{1}{30})(n + \frac{2}{30})\ldots(n + \frac{30}{30})(n+1)}
        {15^{15}10^{10}6^6(n + \frac{1}{15})\ldots(n + \frac{15}{15})(n + \frac{1}{10})\ldots(n + \frac{10}{10})(n + \frac{1}{6})\ldots(n+\frac{6}{6})}.
\]
Now there is a lot of clear cancellation in the fraction, and we see that this is
\[
\frac{30^{30}(n + \frac{1}{30})(n+\frac{7}{30})(n + \frac{11}{30})
        (n + \frac{13}{30})(n + \frac{17}{30})(n + \frac{19}{30})(n + \frac{23}{30})(n + \frac{29}{30})}
    {15^{15}10^{10}6^6 (n + \frac{1}{5})(n + \frac{1}{3})(n + \frac{2}{5})
            (n + \frac{1}{2})(n + \frac{3}{5})(n + \frac{2}{3})(n + \frac{4}{5})(n + 1)}
\]
Which tells us that
\[
    \sum_{n \ge 1} u_n z^n = {}_8F_7\left({\textstyle \frac{1}{30},\frac{7}{30},\frac{11}{30},\frac{13}{30},\frac{17}{30},\frac{19}{30},
            \frac{23}{30},\frac{29}{30};\frac{1}{5},\frac{1}{3},\frac{2}{5},\frac{1}{2},\frac{3}{5},\frac{2}{3},\frac{4}{5}};
            Cz\right),
\]
where
\[
C = \frac{30^{30}}{15^{15}10^{10}6^6}.
\]
\end{exa}

We will need to know that the hypergeometric series attached to a factorial
ratio is essentially unique. We prove this in the next two lemmas.

\begin{lem}\label{misc-lemma-one}
 Suppose that $a_1 \ge a_2 \ge \dots a_K > 0$, $b_1 \ge b_2 \ge \dots, b_L > 0$ and that
\[
 \uab = \uabfrac = 1
\]
for all $n \ge 1$. Then $K=L$ and $\a = \b$.
\end{lem}
\begin{proof}
 The cases $K \ge L$ and $K \le L$ are symmetric, so we may as well assume that
$K \le L$. We will prove the case $K=1$ and then proceed by induction on $K$.

If $K=1$ and $a_1 < b_1$, then it is clear that $u_n \rightarrow 0$ as $n \rightarrow \infty$.
On the other hand, if $a_1 > b_1$, then by Dirichlet's Theorem on primes in
arithmetic progressions, there exists some $m > 1$ such that $a_1m - 1 = p$
is prime. Then $p$ divides the numerator of $u_m$ but not the denominator, so $u_m \ne 1$.
Now, if $a_1 = b_1$, then it is clear that $L = K = 1$.

The case for general $K$ proceeds similarly. We need only show that $a_1 = b_1$, and we are
finished by induction. Again, if $a_1 > b_1$, then there is some $m > 1$ such that
$a_1m - 1 = p$ is prime, and $p$ divides the numerator of $u_m$ but not the denominator.
If, on the other hand, $b_1 > a_1$, we just reverse the argument and find an $m$ and $p$
such that $p$ divides the denominator of $u_m$ but not the numerator.
\end{proof}

\begin{lem}\label{uniqueness-lemma}
 The map
\[
 (\a,\b) \rightarrow u_n(\a,\b;z)
\]
is one-to-one on the set of pairs $(\a,\b)$ such that $a_k \ne  b_l$ for all $k,l$
and $a_1 \ge a_2 \dots \ge a_K$, $b_1 \ge b_2 \ge \ldots \ge b_L$.
\end{lem}
\begin{proof}
 For some $(\a, \b)$ and $(\a', \b')$, we have
\[
 u_n(\a,\b;z) = u_n(\a',\b';z)
\]
if and only if
\[
 u_n(\a,\b) = u_n(\a',\b')
\]
for all $n$. In this case, we can rewrite this as
\[
 u_n(\a,\b)(u_n(\a',\b'))^{-1} = u_n(\a \cup \b', \b \cup \a') = 1.
\]
Now it follows from Lemma \ref{misc-lemma-one} that $\a \cup \b'$
is a permutation of $\b \cup \a'$. Thus $\a'$ is a permutation of $\a$
and $\b'$ is a permutation of $\b$.
\end{proof}
\begin{rek}
 It is also possible to state Lemma \ref{uniqueness-lemma} in an algorithmic
manner. Roughly speaking, given parameters $\alpha = (\alpha_1, \ldots, \alpha_d)$
and $\beta = (\beta_1, \dots, \beta_{d-1})$ that come from a factorial ratio, 
we can form the polynomials $P(x)$ and $Q(x)$ from Definition \ref{P-Q-definition}.
It is then possible to add extra factors to $P(x)$ and $Q(x)$ together to obtain
\[
    \frac{P(x)}{Q(x)} = \frac{(x^{a_1} -1)(x^{a_2} - 1)\cdots(x^{a_K} - 1)}
        {(x^{b_1} - 1)(x^{b_2} - 1)\cdots(x^{b_L} - 1)}
\]
and to recover $\a$ and $\b$. In this manner, if we did not already know about
the $52$ sporadic integer factorial ratio sequences from Vasyunin's work, we could recover them
from the work of Beukers and Heckman \cite{BH-hypergeometric-function} described in
Section \ref{monodromy-section}.
\end{rek}

\printlineno
\printlineno

The main interest in looking at hypergeometric series attached to factorial ratio sequences
comes from the following observation of Rodriguez-Villegas \cite{R-V-factorials}.

\begin{thm}[Rodriguez-Villegas \cite{R-V-factorials}]\label{r-v-algebraic-thm}
Let
\[
    \uab = \uabfrac
\]
with $\sum_{K=1}^K a_k = \sum_{l=1}^L b_l$ and let
\[
    \u(\a,\b; z) = \sum_{n=0}^\infty \uab z^n.
\]
Then $\u(\a,\b;z)$ is algebraic over $\Q(t)$ if and only if $L-K = 1$ and
$\uab \in \Z \text{ for all } n \ge 0.$
\end{thm}

\begin{proof}[Proof of Theorem~{\rm\ref{r-v-algebraic-thm}}, part 1]
We begin by proving that if the generating function is algebraic, then
$u_n$ is in fact integral. In particular, it follows from Lemmas \ref{landau-lemma}
and \ref{k-minus-l-condition} that this will imply that we must have $L-K \ge 1$, which
is, in fact, all that we need from this part of the proof.

A theorem of Eisenstein (see \cite{eisenstein-algebraic-series})
asserts that if $\u(\a;\b;z)$ is algebraic, then
there exists an $N$ such that $\uab \cdot N^n$ in an integer for all $n$. But Lemma
\ref{uab-denominator-lemma} implies that the set of primes occurring in the denominator
of some $\uab$ is either empty or infinite. So, if such an $N$ exists, then
we are able to take $N=1$, which implies that $\uab$ is an integer
for all $n$.
\end{proof}

\printlineno

The remainder of this proof relies on Landau's theorem and the following lemma of
Beukers and Heckman \cite{BH-hypergeometric-function}.
\begin{lem}[Beukers--Heckman \cite{BH-hypergeometric-function}]
Let $\sequence{\alpha}{n}$ and $\sequence{\beta}{n-1}$ be rational numbers
with common denominator $M$. The hypergeometric function
\[
    { }_nF_{n-1}(\sequence{\alpha}{n};\sequence{\beta}{n-1};z)
\]
is algebraic if and only if for all $k$ relatively prime to $M$ the sequences
\[
e(k\alpha_1), \ldots, e(k\alpha_n)
\]
and
\[
e(k\beta_1), \ldots, e(k\beta_{n-1}), 1
\]
interlace on the unit circle.
\end{lem}
\begin{proof}
This follows from \cite[Theorem 4.8]{BH-hypergeometric-function} and the fact that
this function is algebraic if and only if its monodromy group is finite.
\end{proof}

\newcommand{\K}{\kern-.07558875 pt}

For \K the \K case \K of \K the \K hypergeometric \K functions \K that \K are \K generating \K series \K for \K $\uab$
we can make this lemma slightly stronger.


\begin{lem}\label{misc-lemma-two}
\[
\u(\a;\b;z) = \sum_{n=0}^\infty u_n(\a,\b) z^n
\]
is an algebraic function if and only if $\alpha = \alpha(\a,\b)$ and $\beta = \beta(\a,\b)$ interlace
on $[0,1]$.
\end{lem}
\begin{proof}
In our case the $\alpha_i$ and $\beta_j$ are rational numbers in $(0,1]$. Suppose that they have common
denominator $M$. Recall that the numbers $e(\alpha_i)$ are roots of the polynomial $P(\a,\b;x)$, and that
$P$ is the product of cyclotomic polynomials, say $P = \Phi_{m_1}\Phi_{m_2}\cdots\Phi_{m_l}$.
Then for any $(k,M) = 1$ we also have $(k,m_i) = 1$ for all $m_i$. So the map $\alpha \rightarrow \alpha^k$
simply permutes the roots of any $\Phi_{m_i}$. In particular, it permutes the roots of $P$, and hence
permutes the numbers $e(\alpha_i)$.

The exact same argument applies for $\beta$. Thus we have that $\alpha$ and $\beta$ interlace on
$[0,1]$ if and only if $e(\alpha_i)^k$ and $e(\beta_j)^k$ interlace on the unit circle for
all $k$ with $(k,M) = 1$.

In particular, $\u(\a;\b;z)$ is algebraic if and only if $\alpha$ and $\beta$ interlace on
$[0,1]$.
\end{proof}

\printlineno

Combining these lemmas finishes the proof of \ref{r-v-algebraic-thm}.

\begin{proof}[Proof of Theorem~{\rm\ref{r-v-algebraic-thm}}, part 2]
Suppose that $\u(\a,\b;z)$ is algebraic. Then we know that $L-K \ge 1$,
from the first part of the proof. Note that the number of copies of the
number $1$ in the set $\beta$ is $L-K$. However, if $\alpha$ and $\beta$ are to
interlace, no values can be repeated, so we must have $L-K = 1$.

Now, if $L-K = 1$, then from Lemma \ref{equivalence-lemma} we know that $\uab$ is
integral if and only if $\alpha$ and $\beta$ interlace. From Lemma \ref{misc-lemma-two}
we know that this is equivalent to $\u(\a,\b;z)$ being algebraic.
\end{proof}

\section{A Classification of Integral Factorial Ratios}\label{classification-section}
\subsection{Monodromy for Hypergeometric Functions ${ }_nF_{n-1}$}\label{monodromy-section}
This section is an application of the work of Beukers and Heckman \cite{BH-hypergeometric-function},
so we begin by restating a few necessary theorems and definitions.

\newcommand{\w}{\mathbf{w}}
\newcommand{\z}{\mathbf{z}}

\begin{defi}[Hypergeometric Groups]
Let $w_1, \ldots, w_n$ and $z_1, \ldots z_n$ be complex numbers with $w_i \ne z_j$ for all $i$ and
$j$. The hypergeometric group $H(\w,\z)$ with numerator parameters $w_1, \ldots, w_n$ and
denominator parameters $z_1, \ldots, z_n$ is a subgroup of $GL_n(\C)$ generated by elements
\[
    h_0, h_1, \text{ and, } h_\infty
\]
such that
\[
    h_0h_1h_\infty = 1
\]
and
\[
    \det(t - h_\infty) = \prod_{j=1}^n (t - w_j)
\]
\[
    \det(t - h_0^{-1}) = \prod_{j=1}^n (t - z_j)
\]
and such that $h_1 - 1$ has rank $1$.
\end{defi}

\printlineno

Hypergeometric groups are precisely those groups which occur as monodromy groups for
hypergeometric functions. Specifically, we have the following.

\begin{prop}
The monodromy group for the hypergeometric function 
\[
    { }_nF_{n-1}(\alpha_1, \ldots, \alpha_n; \beta_1, \ldots, \beta_{n-1}; z)
\]
is a hypergeometric group with numerator parameters
\[
    e(\alpha_1), e(\alpha_2), \ldots, e(\alpha_n)
\]
and denominator parameters
\[
    e(\beta_1), e(\beta_2), \ldots, e(\beta_{n-1}), 1.
\]
\end{prop}

\begin{proof}
This is \cite[Proposition 3.2]{BH-hypergeometric-function}.
\end{proof}

In categorizing hypergeometric groups it is useful to consider the following special subgroup.
\begin{defi}
The subgroup $H_r(\w,\z)$ of $H(\w,\z)$ generated by $h_\infty^k h_1 h_\infty^{-k}$ for $k \in \Z$
is called the reflection subgroup of $H(\w,\z)$.
\end{defi}

The existence of the following two theorems explains some of the usefulness of considering the
reflection subgroup of a hypergeometric group.
\begin{thm}\label{Hr-reducible-theorem}
The reflection subgroup $H_r(\w,\z)$ of $H(\w,\z)$ acts reducibly on $\C^n$ if and only if there
exists a root of unity $\zeta \ne 1$ such that multiplication by $\zeta$ permutes both the
elements of $\w$ and the elements of $\z$. Moreover, if $H_r(\w,\z)$ is reducible, then
$H(\w,\z)$ is imprimitive.
\end{thm}
\begin{proof}
This is \cite[Theorem 5.3]{BH-hypergeometric-function}.
\end{proof}
\begin{thm}\label{H-primitive-theorem}
Suppose that $H_r(\w,\z)$ is irreducible. Then $H$ is imprimitive if and only if there exist
$p,q \in N$ with $p+q = n$ and $(p,q) = 1$, and $A,B,C \in \C^*$ such that $A^n = B^pC^q$ and
such that
\[
\{w_1, \ldots, w_n\} = \{A, A\zeta_n, A\zeta_n^2, \ldots, A\zeta_n^{n-1}\}
\]
and
\[
\{z_1, \ldots, z_n\} = \{B, B\zeta_p, B\zeta_p^2, \ldots, B\zeta_p^{p-1},
    C, C\zeta_q, C\zeta_q^2, \ldots, C\zeta_q^{q-1}\}
\]
where $\zeta_n = e(1/n)$, or with the same sets of equalities with $\w$ and $\z$ reversed.
\end{thm}
\begin{proof}
 This is \cite[Theorem 5.8]{BH-hypergeometric-function}.
\end{proof}

As defined, hypergeometric groups are subgroups of $GL_n(\C)$. The
following proposition tells us when a hypergeometric group is defined over
$GL_n(R)$ for $R \subset \C$.

\begin{prop}\label{field-of-definition-prop}
 Suppose $w_1, \dots, w_n, z_1, \dots, z_n \in \C^*$ with $w_i \ne z_j$ for all
$i,j$. Let $A_1, \dots A_n, B_1, \dots, B_n$ be defined by
\[
 \prod_{j=1}^n (t - w_j) = t^n + A_1t^{n-1} + \dots + A_n,
\]
and
\[
 \prod_{j=1}^n (t - z_j) = t^n + B_1t^{n-1} + \dots + B_n.
\]
Then relative to a suitable basis, the hypergeometric group $H(\w,\z)$
is defined over the ring $\Z[A_1, \dots, A_n, B_1, \dots, B_n, A_n^{-1}, B_n^{-1}]$.
\end{prop}
\begin{proof}
 This is \cite[Corollary 3.6]{BH-hypergeometric-function}, and follows directly
from a theorem of Levelt \cite[Theorem 1.1]{levelt-thesis}.
\end{proof}

We need to state one more definition, and then we will be ready to state the main
classification theorem of Beukers and Heckman that we are interested in.

\begin{defi}
 A scalar shift of the hypergeometric group $H(\w,\z)$ is a hypergeometric
group $H(d\w,d\z) = H(dw_1,dw_2,\dots,dw_n;dz_1,dz_2,\dots,dz_n)$ for some $d \in C^*$.
\end{defi}

Our main interest in the work of Beukers and Heckman comes from the following
Theorem.

\begin{thm}\label{finite-HG-group-classification}
 Let $n \ge 3$ and let $H(\w,\z) \subset GL_n(\C)$ be a primitive hypergeometric group
whose parameters are roots of unity and generate the field $\Q(\zeta_h)$. Then
$H(\w,\z)$ is finite if and only if, up to a scalar shift, the parameters have the form
$w_1^k, w_2^k, \dots, w_n^k; z_1^k, z_2^k, \dots, z_n^k$, where $\gcd(h,k) = 1$ and
the exponents of either $w_1, \dots, w_n; z_1, \dots, z_n$ or
$z_1, \dots, z_n; w_1, \dots w_n$ are listed in \cite[Table 8.3]{BH-hypergeometric-function}.
\end{thm}
\begin{proof}
 This is \cite[Theorem 7.1]{BH-hypergeometric-function}.
\end{proof}

\subsection{The Classification} From now on we set $L = K+1$. We are interested in ratios
where the parameters
have greatest common divisor $1$. The following lemma shows that this condition translates
nicely into the reflection group of the monodromy group being irreducible.

\begin{lem}
Let $u_n = \uab$ and $\u = \u(\a,\b;z)$. Let $H(\u)$ be the hypergeometric group associated to
$\u$ and let $H_r(\u)$ be the reflection subgroup of $H(\u)$. Suppose that $u_n$ is an integer
for all $n$. Then $H_r(\u)$ acts reducibly on $\C$ if and only if 
\[\gcd(a_1, a_2, \ldots, a_K, b_1, b_2, \ldots, b_{K+1}) > 1.\] Moreover, if $H_r(\u)$ is
reducible, then $H(\u)$ is imprimitive.
\end{lem}

\begin{proof}
Let $P = P(\a,\b;x)$ and $Q = Q(\a,\b;x)$. Then we have
\[
    \frac{P}{Q} = \frac{(x - e(\alpha_1))\cdots(x - e(\alpha_d))}{(x - e(\beta_1))\cdots(x - e(\beta_d))}
        = \frac{(x^{a_1} - 1)\cdots(x^{a_K} - 1)}{(x^{b_1} - 1)\cdots(x^{b_{K+1}} - 1)},
\]
and $H(\u)$ is a hypergeometric group with numerator parameters 
\[
e(\alpha_1), e(\alpha_2), \ldots, e(\alpha_d)
\]
and
denominator parameters
\[
e(\beta_1), e(\beta_2), \ldots, e(\beta_d).
\]
From Theorem \ref{Hr-reducible-theorem} we know that
$H_r(\u)$ acts reducibly on $\C$ if and only if there exists some $\gamma \not\equiv 0 \mod 1$ such
that
\[
    \{e(\alpha_1), e(\alpha_2), \ldots, e(\alpha_d)\} = 
        \{e(\alpha_1 + \gamma), e(\alpha_2 + \gamma), \ldots, e(\alpha_d + \gamma)\}
\]
and
\[
    \{e(\beta_1), e(\beta_2), \ldots, e(\beta_d)\} = 
        \{e(\beta_1 + \gamma), e(\beta_2 + \gamma), \ldots, e(\beta_d + \gamma)\}.
\]
Such a permutation by multiplication by $e(\gamma)$ is a rotation of the unit circle, and necessarily
lifts to permutations of the roots of $(x^{a_1} - 1)(x^{a_2} - 1)\cdots(x^{a_K} - 1)$ and the
roots of $(x^{b_1} - 1)(x^{b_2} - 1)\cdots(x^{b_{K+1}} - 1)$. This permutation takes the root
$e(1)$ to the root $e(\gamma)$, so $e(\gamma)$ must be a root of $(x^{a_i} - 1)$ and a root
of $(x^{b_j} - 1)$ for all $i$ and $j$. If $e(\gamma)$ is a primitive $M$th root of unity,
then this means that $M | a_i$ and $M | b_j$ for all $i$ and $j$. Thus, such a nonzero
$\gamma$ exists if and only if the $\gcd$ of all of the terms greater than $1$.
\end{proof}

\printlineno

Vasyunin noticed that the step functions corresponding to
\begin{equation}\label{five-term-family-one}
u_n = \frac{(2an)!(bn)!}{(an)!(2bn)![(a-b)n]!} \text{ with } b < a
\end{equation}
and
\begin{equation}\label{five-term-family-two}
u_n = \frac{(2an)!(2bn)!}{(an)!(bn)![(a+b)n]!}
\end{equation}
are nonnegative. Thus for both of these families,
$u_n$ is an integer for all $n$. It turns out that when $a$ and $b$ are not both
odd these two infinite families give
exactly those with factorial ratios with $\gcd 1 $ for which the hypergeometric group is imprimitive.
On the other hand, when $a$ and $b$ are both odd, these come from
scalar shifts of the hypergeometric groups associated to binomial coefficients.

\begin{lem}\label{reflection-group-irreducible-lemma}
Let $u_n = \uab$ and $\u = \u(\a,\b;z)$. Let $H(\u)$ be the hypergeometric group associated to
$\u$ and let $H_r(\u)$ be the reflection subgroup of $H(\u)$. Suppose that $u_n$ is an integer
for all $n$ and that $H_r(\u)$ is irreducible. Then $H(\u)$ is imprimitive if and only
if $u_n$ is of the form \eqref{five-term-family-one} or \eqref{five-term-family-two} with
$a$ and $b$ not both odd.
\end{lem}

\begin{proof}
Again let $P = P(\a,\b;x)$ and $Q = Q(\a,\b;x)$. Suppose that $u_n$ is the of form
\eqref{five-term-family-one}. Then we have
\[
    \frac{P}{Q} = \frac{(x^{2a} - 1)(x^b - 1)}{(x^{2b} - 1)(x^a - 1)(x^{a-b} - 1)}
        = \frac{(x^a + 1)}{(x^b + 1)(x^{a-b} - 1)},
\]
which is in lowest terms if $a$ and $b$ are not both odd.
Then the numerator parameters for $H(\u)$ are
\[
    A, A\zeta_a^1, A\zeta_a^2, \ldots, A\zeta_a^{a-1}
\]
and the denominator parameters are
\[
B, B\zeta_b, B\zeta_b^2, \ldots, B\zeta_b^{b-1},
    \zeta_{a-b}, \zeta_{a-b}^2, \ldots, \zeta_{a-b}^{a-b-1}, 1
\]
where $A = \zeta_{2a}$ and $B = \zeta_{2b}$ satisfy $A^a = B^b = -1$, so these parameters satisfy
the condition of Theorem \ref{H-primitive-theorem}, and $H(\u)$ is imprimitive.

Similarly, if $u_n$ is of the form \eqref{five-term-family-two} then we have
\[
    \frac{P}{Q} = \frac{(x^{2a} - 1)(x^{2b} - 1)}{(x^a - 1)(x^b - 1)(x^{a+b} - 1)}
        = \frac{(x^a + 1)(x^b + 1)}{x^{a+b} - 1},
\]
again in lowest terms if $a$ and $b$ are not both odd,
and so $H(\u)$ is a hypergeometric group with numerator parameters
\[
    A, A\zeta_a^1, A\zeta_a^2, \ldots, A\zeta_a^{a-1},
    B, B\zeta_b, B\zeta_b^2, \ldots, B\zeta_b^{b-1}
\]
and denominator parameters
\[
    \zeta_{a+b},\zeta_{a+b}^2, \ldots, \zeta_{a+b}^{a+b}
\]
where $A = \zeta_{2a}$ and $B = \zeta_{2b}$. $A$ and $B$ satisfy $A^a = B^b$, so
$H(\u)$ again satisfies the conditions of the theorem.

To see the converse, suppose that the numerator parameters of $H(u)$ are of the form
\[
A, A\zeta_a, A\zeta_a^2, \ldots, A\zeta_a^{a-1}.
\]
These parameters must have the property
that, if they contain one primitive $M$th root of unity for some $M$, then they
contain all of them. Thus, by symmetry considerations, we find that the only possibility
is that $A = \zeta_{2a}$. Similarly, for the denominator parameters we find that
either $B$ or $C$ is $1$, and so without loss of generality, we assume that $C = 1$ and
find that $B = \zeta_{2b}$ is the only possibility. Indeed, whenever $(a,b) = 1$ and $a$ and
$b$ are not both odd, this does
work and gives $u_n$ of the form \eqref{five-term-family-one}.

A similar argument works for the second case.
\end{proof}

We now examine the case where both $a$ and $b$ are odd.

\begin{lem}\label{binomial-shift-lemma}
 Let $u_n = \uab$ and $\u = \u(\a,\b;z)$. Let $H(\u)$ be the hypergeometric
group associated to $\u$. 
\begin{itemize}
 \item[(i)] If $u_n$ is of the form \eqref{five-term-family-one} with $a$ and $b$ both odd, then $H(\u)$
is a scalar shift by $-1 = e(1/2)$ of $H(\u')$, where
\[
 u_n' = {an \choose bn}.
\]
\item[(ii)] If $u_n$ is of the form \eqref{five-term-family-two} with $a$ and $b$ both odd, then
$H(\u)$ is obtained by taking a scalar shift of $H(\u')$ and reversing the
numerator and denominator parameters, where
\[
 u_n' = {(a+b)n \choose an}
\]
\end{itemize}
\end{lem}
\begin{proof}
 (i) Suppose that $u_n$ is of the form \eqref{five-term-family-one}. Then, as in the previous
lemma, we have for $P(x) = P(\a,\b;x)$ and $Q(x) = Q(\a,\b;x)$,
\[
 \frac{P(x)}{Q(x)} = \frac{(x^a + 1)}{(x^b + 1)(x^{a-b} - 1)}.
\]
This is not in lowest terms, but it is in a convenient form for computing the scalar shift of
$H(\u)$.  The scalar shift corresponds to multiplying each root of $P$ and $Q$ by $-1$, in which
case we obtain polynomials $P^*$ and $Q^*$ with
\[
 \frac{P^*(x)}{Q^*(x)} = \frac{(x^a - 1)}{(x^b - 1)(x^{a-b} - 1)},
\]
which very clearly come from $u_n' = {an \choose bn}$.

(ii) If $u_n$ is of the form \eqref{five-term-family-two}, we proceed similarly, except that
this time we find that
\[
 \frac{P^*(x)}{Q^*(x)} = \frac{(x^a - 1)(x^b - 1)}{(x^{a+b} - 1)},
\]
which we can see comes from $u_n' = {(a+b)n \choose an}$, with the numerator
and denominator parameters reversed.
\end{proof}


We are now ready to prove Theorem \ref{main-theorem}

\begin{proof}[Proof of Theorem~{\rm\ref{main-theorem}}]
Lemma \ref{reflection-group-irreducible-lemma} classifies all of those integral
factorial ratios whose associated hypergeometric group is imprimitive, so it remains to
classify those with associated primitive hypergeometric groups. Beukers and Heckman have
categorized all finite primitive hypergeometric groups,
so we can examine \cite[Table 8.3]{BH-hypergeometric-function} to find all integral factorial
ratios with associated primitive hypergeometric groups.

Specifically, it follows from directly from Proposition \ref{field-of-definition-prop}
that a hypergeometric group associated to
a factorial ratio must be defined over $\Q$.
There are $27$ entries in \cite[Table 8.3]{BH-hypergeometric-function}
that are defined over $\Q$. So, if $H$ is a primitive hypergeometric group which comes from
an integral factorial ratio sequence, then $H$ is either one of these entries or a scalar
shift of one of these entries, possibly with the numerator and denominator parameters switched.
Moreover, the only scalar shift that we need consider is a scalar shift
by $-1$, as a scalar shift by any other root of unity will change the field of
definition.

This yields an infinite family of parameter sets from Line $1$ of \cite[Table 8.3]{BH-hypergeometric-function},
and $104$ ``sporadic'' possibilities. However, the denominator parameters of
a hypergeometric group coming from a factorial ratio must contain a $1$, so we really
only have $52$ possibilities. It is easy to computationally verify that each of the
$52$ sporadic step functions found by Vasyunin corresponds to one of these
$52$ hypergeometric group parameter sets, but this step is not even
necessary, since Lemma \ref{uniqueness-lemma} assures us that different
factorial ratios must have different hypergeometric group parameter sets.

It is easily seen that Line 1 of \cite[Table 8.3]{BH-hypergeometric-function} corresponds to the infinite family
of binomial coefficient sequences
\[
	\frac{[(a+b)n]!}{(an)!(bn)!} \textrm{ with } \gcd(a,b) = 1,
\]
and, as we have seen in Lemma \ref{binomial-shift-lemma},
scalar shifts of this family by $-1$ yield factorial ratios of the forms
\eqref{five-term-family-one} and \eqref{five-term-family-two} already considered.

Lemma \ref{uniqueness-lemma} assures us that this must be all integer factorial
ratio sequences in which the parameters have greatest common divisor $1$, giving
us the ``only if'' part of the theorem.
\end{proof}

\section{A Listing of all Integral Factorial Ratios with $L-K = 1$}
The following tables contain a listing of all solutions to Problem \ref{uab-problem} with
$\gcd 1$, organized as follows. The second column lists the parameters for $\uab$. The parameter
$d$ of the third column is the dimension of the monodromy group of $\u(\a,\b;z)$, and the
fourth column lists the specific parameters for ${ }_dF_{d-1}$. With the exception of lines
2 and 3 of Table \ref{factorial-ratio-table-two} all the entries have primitive monodromy groups,
and so they have a corresponding entry in \cite[Table 8.3]{BH-hypergeometric-function}.

It should be noted that lines 1, 2, and 3 of Table \ref{factorial-ratio-table-two} only correspond to solutions with $\gcd 1$ when
$\gcd(a,b) = 1$.
\begin{longtable}{Sl|Sc|Sc|Sc|Sr}
\caption{The three infinite families of integral factorial ratio sequences}\label{factorial-ratio-table-two} \\
Line \#& $\uab$& $\phantom{}d\phantom{}$ & ${ }_dF_{d-1}$ parameters & \cite{BH-hypergeometric-function} Line \#\\
\hline\hline&&&\\
\endfirsthead
Line \#& $\uab$& $\phantom{}d\phantom{}$ & ${ }_dF_{d-1}$ parameters & \cite{BH-hypergeometric-function} Line \#\\
\hline\hline&&&\\
\endhead
$1_{a,b}$ & $\begin{array}{Sc} \protect{$\scriptstyle[a+b]$} \\ \protect{$\scriptstyle[a,b]$} \end{array}$ & $\!a+b-1\!$ &
    $\begin{array}{c} \protect{\scriptstyle[\frac{1}{a+b},\frac{2}{a+b},\ldots,\frac{a+b-1}{a+b}]} \\
    \protect{\scriptstyle[\frac{1}{a},\frac{2}{a},\ldots,\frac{a-1}{a},\frac{1}{b},\frac{2}{b},\ldots,\frac{b-1}{b}]}
    \end{array}$ & $1$ \\
$2_{a,b}$ & $\begin{array}{Sc} \protect{$\scriptstyle [2a,b]$} \\ \protect{$\scriptstyle [a,2b,a-b]$} \end{array}$ & $a$ &
    $\begin{array}{Sc}\protect{$\scriptstyle[\frac{1}{2a},\frac{3}{2a},\ldots\frac{2a-1}{2a}]$} \\
        \protect{$\scriptstyle[\!\frac{1}{2b},\frac{3}{2b},\ldots\frac{2b-1}{2b},\frac{1}{a-b},\frac{2}{a-b},\ldots,
        \frac{a-b-1}{a-b}]\!$}\end{array}$ & None\\
$3_{a,b}$ & $\begin{array}{Sc} \protect{$\scriptstyle[2a,2b]$} \\ \protect{$\scriptstyle[a,b,a+b]$} \end{array}$ & $a+b$ &
    $\begin{array}{Sc}\protect{$\scriptstyle[\frac{1}{2a},\frac{3}{2a},\ldots\frac{2a-1}{2a},\frac{1}{2b},\frac{3}{2b},\ldots\frac{2b-1}{2b}]$} \\
        \protect{$\scriptstyle[\frac{1}{a+b},\frac{2}{a+b},\ldots,
        \frac{a+b-1}{a+b}]$}\end{array}$ & None\\
\end{longtable}
\begin{longtable}{Sl|Sc|Sc|Sc|Sr}
\caption{The $52$ sporadic integral factorial ratio sequences}\label{factorial-ratio-table} \\
Line \#& $\uab$& $\phantom{M}d\phantom{M}$ & \phantom{MM}${ }_dF_{d-1}$ parameters\phantom{MM}& \cite{BH-hypergeometric-function} Line \#\\
\hline\hline&&&\\
\endfirsthead
Line \#& $\uab$& $\phantom{M}d\phantom{M}$ & \phantom{MM}${ }_dF_{d-1}$ parameters\phantom{MM}& \cite{BH-hypergeometric-function} Line \#\\
\hline\hline&&&\\
\endhead

            $1$ &
            $\begin{array}{Sc}\protect{$\scriptstyle [12, 1]$}\\ \protect{$\scriptstyle [6, 4, 3]$} \end{array}$ &
            $4$ &
            $\begin{array}{Sc}\protect{$\scriptstyle \left[\frac{1}{12}, 
 \frac{5}{12}, 
 \frac{7}{12}, 
 \frac{11}{12}\right]$} \\ \protect{$\scriptstyle \left[\frac{1}{3}, 
 \frac{1}{2}, 
 \frac{2}{3}\right]$} \end{array}$ &
            $37$
            \\

            $2$ &
            $\begin{array}{Sc}\protect{$\scriptstyle [12, 3, 2]$}\\ \protect{$\scriptstyle [6, 6, 4, 1]$} \end{array}$ &
            $4$ &
            $\begin{array}{Sc}\protect{$\scriptstyle \left[\frac{1}{12}, 
 \frac{5}{12}, 
 \frac{7}{12}, 
 \frac{11}{12}\right]$} \\ \protect{$\scriptstyle \left[\frac{1}{6}, 
 \frac{1}{2}, 
 \frac{5}{6}\right]$} \end{array}$ &
            $37$
            \\

            $3$ &
            $\begin{array}{Sc}\protect{$\scriptstyle [12, 1]$}\\ \protect{$\scriptstyle [8, 3, 2]$} \end{array}$ &
            $6$ &
            $\begin{array}{Sc}\protect{$\scriptstyle \left[\frac{1}{12}, 
 \frac{1}{6}, 
 \frac{5}{12}, 
 \frac{7}{12}, 
 \frac{5}{6}, 
 \frac{11}{12}\right]$} \\ \protect{$\scriptstyle \left[\frac{1}{8}, 
 \frac{3}{8}, 
 \frac{1}{2}, 
 \frac{5}{8}, 
 \frac{7}{8}\right]$} \end{array}$ &
            $45$
            \\

            $4$ &
            $\begin{array}{Sc}\protect{$\scriptstyle [12, 3]$}\\ \protect{$\scriptstyle [8, 6, 1]$} \end{array}$ &
            $6$ &
            $\begin{array}{Sc}\protect{$\scriptstyle \left[\frac{1}{12}, 
 \frac{1}{3}, 
 \frac{5}{12}, 
 \frac{7}{12}, 
 \frac{2}{3}, 
 \frac{11}{12}\right]$} \\ \protect{$\scriptstyle \left[\frac{1}{8}, 
 \frac{3}{8}, 
 \frac{1}{2}, 
 \frac{5}{8}, 
 \frac{7}{8}\right]$} \end{array}$ &
            $45$
            \\

            $5$ &
            $\begin{array}{Sc}\protect{$\scriptstyle [12, 3]$}\\ \protect{$\scriptstyle [6, 5, 4]$} \end{array}$ &
            $6$ &
            $\begin{array}{Sc}\protect{$\scriptstyle \left[\frac{1}{12}, 
 \frac{1}{3}, 
 \frac{5}{12}, 
 \frac{7}{12}, 
 \frac{2}{3}, 
 \frac{11}{12}\right]$} \\ \protect{$\scriptstyle \left[\frac{1}{5}, 
 \frac{2}{5}, 
 \frac{1}{2}, 
 \frac{3}{5}, 
 \frac{4}{5}\right]$} \end{array}$ &
            $46$
            \\

            $6$ &
            $\begin{array}{Sc}\protect{$\scriptstyle [12, 5]$}\\ \protect{$\scriptstyle [10, 4, 3]$} \end{array}$ &
            $6$ &
            $\begin{array}{Sc}\protect{$\scriptstyle \left[\frac{1}{12}, 
 \frac{1}{6}, 
 \frac{5}{12}, 
 \frac{7}{12}, 
 \frac{5}{6}, 
 \frac{11}{12}\right]$} \\ \protect{$\scriptstyle \left[\frac{1}{10}, 
 \frac{3}{10}, 
 \frac{1}{2}, 
 \frac{7}{10}, 
 \frac{9}{10}\right]$} \end{array}$ &
            $46$
            \\

            $7$ &
            $\begin{array}{Sc}\protect{$\scriptstyle [18, 1]$}\\ \protect{$\scriptstyle [9, 6, 4]$} \end{array}$ &
            $6$ &
            $\begin{array}{Sc}\protect{$\scriptstyle \left[\frac{1}{18}, 
 \frac{5}{18}, 
 \frac{7}{18}, 
 \frac{11}{18}, 
 \frac{13}{18}, 
 \frac{17}{18}\right]$} \\ \protect{$\scriptstyle \left[\frac{1}{4}, 
 \frac{1}{3}, 
 \frac{1}{2}, 
 \frac{2}{3}, 
 \frac{3}{4}\right]$} \end{array}$ &
            $47$
            \\

            $8$ &
            $\begin{array}{Sc}\protect{$\scriptstyle [9, 2]$}\\ \protect{$\scriptstyle [6, 4, 1]$} \end{array}$ &
            $6$ &
            $\begin{array}{Sc}\protect{$\scriptstyle \left[\frac{1}{9}, 
 \frac{2}{9}, 
 \frac{4}{9}, 
 \frac{5}{9}, 
 \frac{7}{9}, 
 \frac{8}{9}\right]$} \\ \protect{$\scriptstyle \left[\frac{1}{6}, 
 \frac{1}{4}, 
 \frac{1}{2}, 
 \frac{3}{4}, 
 \frac{5}{6}\right]$} \end{array}$ &
            $47$
            \\

            $9$ &
            $\begin{array}{Sc}\protect{$\scriptstyle [9, 4]$}\\ \protect{$\scriptstyle [8, 3, 2]$} \end{array}$ &
            $6$ &
            $\begin{array}{Sc}\protect{$\scriptstyle \left[\frac{1}{9}, 
 \frac{2}{9}, 
 \frac{4}{9}, 
 \frac{5}{9}, 
 \frac{7}{9}, 
 \frac{8}{9}\right]$} \\ \protect{$\scriptstyle \left[\frac{1}{8}, 
 \frac{3}{8}, 
 \frac{1}{2}, 
 \frac{5}{8}, 
 \frac{7}{8}\right]$} \end{array}$ &
            $48$
            \\

            $10$ &
            $\begin{array}{Sc}\protect{$\scriptstyle [18, 4, 3]$}\\ \protect{$\scriptstyle [9, 8, 6, 2]$} \end{array}$ &
            $6$ &
            $\begin{array}{Sc}\protect{$\scriptstyle \left[\frac{1}{18}, 
 \frac{5}{18}, 
 \frac{7}{18}, 
 \frac{11}{18}, 
 \frac{13}{18}, 
 \frac{17}{18}\right]$} \\ \protect{$\scriptstyle \left[\frac{1}{8}, 
 \frac{3}{8}, 
 \frac{1}{2}, 
 \frac{5}{8}, 
 \frac{7}{8}\right]$} \end{array}$ &
            $48$
            \\

            $11$ &
            $\begin{array}{Sc}\protect{$\scriptstyle [9, 1]$}\\ \protect{$\scriptstyle [5, 3, 2]$} \end{array}$ &
            $6$ &
            $\begin{array}{Sc}\protect{$\scriptstyle \left[\frac{1}{9}, 
 \frac{2}{9}, 
 \frac{4}{9}, 
 \frac{5}{9}, 
 \frac{7}{9}, 
 \frac{8}{9}\right]$} \\ \protect{$\scriptstyle \left[\frac{1}{5}, 
 \frac{2}{5}, 
 \frac{1}{2}, 
 \frac{3}{5}, 
 \frac{4}{5}\right]$} \end{array}$ &
            $49$
            \\

            $12$ &
            $\begin{array}{Sc}\protect{$\scriptstyle [18, 5, 3]$}\\ \protect{$\scriptstyle [10, 9, 6, 1]$} \end{array}$ &
            $6$ &
            $\begin{array}{Sc}\protect{$\scriptstyle \left[\frac{1}{18}, 
 \frac{5}{18}, 
 \frac{7}{18}, 
 \frac{11}{18}, 
 \frac{13}{18}, 
 \frac{17}{18}\right]$} \\ \protect{$\scriptstyle \left[\frac{1}{10}, 
 \frac{3}{10}, 
 \frac{1}{2}, 
 \frac{7}{10}, 
 \frac{9}{10}\right]$} \end{array}$ &
            $49$
            \\

            $13$ &
            $\begin{array}{Sc}\protect{$\scriptstyle [18, 4]$}\\ \protect{$\scriptstyle [12, 9, 1]$} \end{array}$ &
            $7$ &
            $\begin{array}{Sc}\protect{$\scriptstyle \left[\frac{1}{18}, 
\frac{5}{18}, 
\frac{7}{18}, 
\frac{1}{2}, 
\frac{11}{18}, 
\frac{13}{18}, 
\frac{17}{18}\right]$} \\ \protect{$\scriptstyle \left[\frac{1}{12}, 
 \frac{1}{3}, 
 \frac{5}{12}, 
 \frac{7}{12}, 
 \frac{2}{3}, 
 \frac{11}{12}\right]$} \end{array}$ &
            $58$
            \\

            $14$ &
            $\begin{array}{Sc}\protect{$\scriptstyle [12, 2]$}\\ \protect{$\scriptstyle [9, 4, 1]$} \end{array}$ &
            $7$ &
            $\begin{array}{Sc}\protect{$\scriptstyle \left[\frac{1}{12}, 
\frac{1}{6}, 
\frac{5}{12}, 
\frac{1}{2}, 
\frac{7}{12}, 
\frac{5}{6}, 
\frac{11}{12}\right]$} \\ \protect{$\scriptstyle \left[\frac{1}{9}, 
 \frac{2}{9}, 
 \frac{4}{9}, 
 \frac{5}{9}, 
 \frac{7}{9}, 
 \frac{8}{9}\right]$} \end{array}$ &
            $58$
            \\

            $15$ &
            $\begin{array}{Sc}\protect{$\scriptstyle [18, 2]$}\\ \protect{$\scriptstyle [9, 6, 5]$} \end{array}$ &
            $7$ &
            $\begin{array}{Sc}\protect{$\scriptstyle \left[\frac{1}{18}, 
\frac{5}{18}, 
\frac{7}{18}, 
\frac{1}{2}, 
\frac{11}{18}, 
\frac{13}{18}, 
\frac{17}{18}\right]$} \\ \protect{$\scriptstyle \left[\frac{1}{5}, 
 \frac{1}{3}, 
 \frac{2}{5}, 
 \frac{3}{5}, 
 \frac{2}{3}, 
 \frac{4}{5}\right]$} \end{array}$ &
            $59$
            \\

            $16$ &
            $\begin{array}{Sc}\protect{$\scriptstyle [10, 6]$}\\ \protect{$\scriptstyle [9, 5, 2]$} \end{array}$ &
            $7$ &
            $\begin{array}{Sc}\protect{$\scriptstyle \left[\frac{1}{10}, 
\frac{1}{6}, 
\frac{3}{10}, 
\frac{1}{2}, 
\frac{7}{10}, 
\frac{5}{6}, 
\frac{9}{10}\right]$} \\ \protect{$\scriptstyle \left[\frac{1}{9}, 
 \frac{2}{9}, 
 \frac{4}{9}, 
 \frac{5}{9}, 
 \frac{7}{9}, 
 \frac{8}{9}\right]$} \end{array}$ &
            $59$
            \\

            $17$ &
            $\begin{array}{Sc}\protect{$\scriptstyle [14, 3]$}\\ \protect{$\scriptstyle [9, 7, 1]$} \end{array}$ &
            $7$ &
            $\begin{array}{Sc}\protect{$\scriptstyle \left[\frac{1}{14}, 
\frac{3}{14}, 
\frac{5}{14}, 
\frac{1}{2}, 
\frac{9}{14}, 
\frac{11}{14}, 
\frac{13}{14}\right]$} \\ \protect{$\scriptstyle \left[\frac{1}{9}, 
 \frac{2}{9}, 
 \frac{4}{9}, 
 \frac{5}{9}, 
 \frac{7}{9}, 
 \frac{8}{9}\right]$} \end{array}$ &
            $60$
            \\

            $18$ &
            $\begin{array}{Sc}\protect{$\scriptstyle [18, 3, 2]$}\\ \protect{$\scriptstyle [9, 7, 6, 1]$} \end{array}$ &
            $7$ &
            $\begin{array}{Sc}\protect{$\scriptstyle \left[\frac{1}{18}, 
\frac{5}{18}, 
\frac{7}{18}, 
\frac{1}{2}, 
\frac{11}{18}, 
\frac{13}{18}, 
\frac{17}{18}\right]$} \\ \protect{$\scriptstyle \left[\frac{1}{7}, 
 \frac{2}{7}, 
 \frac{3}{7}, 
 \frac{4}{7}, 
 \frac{5}{7}, 
 \frac{6}{7}\right]$} \end{array}$ &
            $60$
            \\

            $19$ &
            $\begin{array}{Sc}\protect{$\scriptstyle [12, 2]$}\\ \protect{$\scriptstyle [7, 4, 3]$} \end{array}$ &
            $7$ &
            $\begin{array}{Sc}\protect{$\scriptstyle \left[\frac{1}{12}, 
\frac{1}{6}, 
\frac{5}{12}, 
\frac{1}{2}, 
\frac{7}{12}, 
\frac{5}{6}, 
\frac{11}{12}\right]$} \\ \protect{$\scriptstyle \left[\frac{1}{7}, 
 \frac{2}{7}, 
 \frac{3}{7}, 
 \frac{4}{7}, 
 \frac{5}{7}, 
 \frac{6}{7}\right]$} \end{array}$ &
            $61$
            \\

            $20$ &
            $\begin{array}{Sc}\protect{$\scriptstyle [14, 6, 4]$}\\ \protect{$\scriptstyle [12, 7, 3, 2]$} \end{array}$ &
            $7$ &
            $\begin{array}{Sc}\protect{$\scriptstyle \left[\frac{1}{14}, 
\frac{3}{14}, 
\frac{5}{14}, 
\frac{1}{2}, 
\frac{9}{14}, 
\frac{11}{14}, 
\frac{13}{14}\right]$} \\ \protect{$\scriptstyle \left[\frac{1}{12}, 
 \frac{1}{3}, 
 \frac{5}{12}, 
 \frac{7}{12}, 
 \frac{2}{3}, 
 \frac{11}{12}\right]$} \end{array}$ &
            $61$
            \\

            $21$ &
            $\begin{array}{Sc}\protect{$\scriptstyle [14, 1]$}\\ \protect{$\scriptstyle [7, 5, 3]$} \end{array}$ &
            $7$ &
            $\begin{array}{Sc}\protect{$\scriptstyle \left[\frac{1}{14}, 
\frac{3}{14}, 
\frac{5}{14}, 
\frac{1}{2}, 
\frac{9}{14}, 
\frac{11}{14}, 
\frac{13}{14}\right]$} \\ \protect{$\scriptstyle \left[\frac{1}{5}, 
 \frac{1}{3}, 
 \frac{2}{5}, 
 \frac{3}{5}, 
 \frac{2}{3}, 
 \frac{4}{5}\right]$} \end{array}$ &
            $62$
            \\

            $22$ &
            $\begin{array}{Sc}\protect{$\scriptstyle [10, 6, 1]$}\\ \protect{$\scriptstyle [7, 5, 3, 2]$} \end{array}$ &
            $7$ &
            $\begin{array}{Sc}\protect{$\scriptstyle \left[\frac{1}{10}, 
\frac{1}{6}, 
\frac{3}{10}, 
\frac{1}{2}, 
\frac{7}{10}, 
\frac{5}{6}, 
\frac{9}{10}\right]$} \\ \protect{$\scriptstyle \left[\frac{1}{7}, 
 \frac{2}{7}, 
 \frac{3}{7}, 
 \frac{4}{7}, 
 \frac{5}{7}, 
 \frac{6}{7}\right]$} \end{array}$ &
            $62$
            \\

            $23$ &
            $\begin{array}{Sc}\protect{$\scriptstyle [15, 1]$}\\ \protect{$\scriptstyle [9, 5, 2]$} \end{array}$ &
            $8$ &
            $\begin{array}{Sc}\protect{$\scriptstyle \left[\frac{1}{15}, 
\frac{2}{15}, 
\frac{4}{15}, 
\frac{7}{15}, 
\frac{8}{15}, 
\frac{11}{15}, 
\frac{13}{15}, 
\frac{14}{15}\right]$} \\ \protect{$\scriptstyle \left[\frac{1}{9}, 
 \frac{2}{9}, 
 \frac{4}{9}, 
 \frac{1}{2}, 
 \frac{5}{9}, 
 \frac{7}{9}, 
 \frac{8}{9}\right]$} \end{array}$ &
            $63$
            \\

            $24$ &
            $\begin{array}{Sc}\protect{$\scriptstyle [30, 9, 5]$}\\ \protect{$\scriptstyle [18, 15, 10, 1]$} \end{array}$ &
            $8$ &
            $\begin{array}{Sc}\protect{$\scriptstyle \left[\frac{1}{30}, 
\frac{7}{30}, 
\frac{11}{30}, 
\frac{13}{30}, 
\frac{17}{30}, 
\frac{19}{30}, 
\frac{23}{30}, 
\frac{29}{30}\right]$} \\ \protect{$\scriptstyle \left[\frac{1}{18}, 
\frac{5}{18}, 
\frac{7}{18}, 
\frac{1}{2}, 
\frac{11}{18}, 
\frac{13}{18}, 
\frac{17}{18}\right]$} \end{array}$ &
            $63$
            \\

            $25$ &
            $\begin{array}{Sc}\protect{$\scriptstyle [15, 4]$}\\ \protect{$\scriptstyle [12, 5, 2]$} \end{array}$ &
            $8$ &
            $\begin{array}{Sc}\protect{$\scriptstyle \left[\frac{1}{15}, 
\frac{2}{15}, 
\frac{4}{15}, 
\frac{7}{15}, 
\frac{8}{15}, 
\frac{11}{15}, 
\frac{13}{15}, 
\frac{14}{15}\right]$} \\ \protect{$\scriptstyle \left[\frac{1}{12}, 
\frac{1}{6}, 
\frac{5}{12}, 
\frac{1}{2}, 
\frac{7}{12}, 
\frac{5}{6}, 
\frac{11}{12}\right]$} \end{array}$ &
            $64$
            \\

            $26$ &
            $\begin{array}{Sc}\protect{$\scriptstyle [30, 5, 4]$}\\ \protect{$\scriptstyle [15, 12, 10, 2]$} \end{array}$ &
            $8$ &
            $\begin{array}{Sc}\protect{$\scriptstyle \left[\frac{1}{30}, 
\frac{7}{30}, 
\frac{11}{30}, 
\frac{13}{30}, 
\frac{17}{30}, 
\frac{19}{30}, 
\frac{23}{30}, 
\frac{29}{30}\right]$} \\ \protect{$\scriptstyle \left[\frac{1}{12}, 
\frac{1}{3}, 
\frac{5}{12}, 
\frac{1}{2}, 
\frac{7}{12}, 
\frac{2}{3}, 
\frac{11}{12}\right]$} \end{array}$ &
            $64$
            \\

            $27$ &
            $\begin{array}{Sc}\protect{$\scriptstyle [15, 4]$}\\ \protect{$\scriptstyle [8, 6, 5]$} \end{array}$ &
            $8$ &
            $\begin{array}{Sc}\protect{$\scriptstyle \left[\frac{1}{15}, 
\frac{2}{15}, 
\frac{4}{15}, 
\frac{7}{15}, 
\frac{8}{15}, 
\frac{11}{15}, 
\frac{13}{15}, 
\frac{14}{15}\right]$} \\ \protect{$\scriptstyle \left[\frac{1}{8}, 
 \frac{1}{6}, 
 \frac{3}{8}, 
 \frac{1}{2}, 
 \frac{5}{8}, 
 \frac{5}{6}, 
 \frac{7}{8}\right]$} \end{array}$ &
            $65$
            \\

            $28$ &
            $\begin{array}{Sc}\protect{$\scriptstyle [30, 5, 4]$}\\ \protect{$\scriptstyle [15, 10, 8, 6]$} \end{array}$ &
            $8$ &
            $\begin{array}{Sc}\protect{$\scriptstyle \left[\frac{1}{30}, 
\frac{7}{30}, 
\frac{11}{30}, 
\frac{13}{30}, 
\frac{17}{30}, 
\frac{19}{30}, 
\frac{23}{30}, 
\frac{29}{30}\right]$} \\ \protect{$\scriptstyle \left[\frac{1}{8}, 
 \frac{1}{3}, 
 \frac{3}{8}, 
 \frac{1}{2}, 
 \frac{5}{8}, 
 \frac{2}{3}, 
 \frac{7}{8}\right]$} \end{array}$ &
            $65$
            \\

            $29$ &
            $\begin{array}{Sc}\protect{$\scriptstyle [15, 2]$}\\ \protect{$\scriptstyle [10, 4, 3]$} \end{array}$ &
            $8$ &
            $\begin{array}{Sc}\protect{$\scriptstyle \left[\frac{1}{15}, 
\frac{2}{15}, 
\frac{4}{15}, 
\frac{7}{15}, 
\frac{8}{15}, 
\frac{11}{15}, 
\frac{13}{15}, 
\frac{14}{15}\right]$} \\ \protect{$\scriptstyle \left[\frac{1}{10}, 
\frac{1}{4}, 
\frac{3}{10}, 
\frac{1}{2}, 
\frac{7}{10}, 
\frac{3}{4}, 
\frac{9}{10}\right]$} \end{array}$ &
            $66$
            \\

            $30$ &
            $\begin{array}{Sc}\protect{$\scriptstyle [30, 3, 2]$}\\ \protect{$\scriptstyle [15, 10, 6, 4]$} \end{array}$ &
            $8$ &
            $\begin{array}{Sc}\protect{$\scriptstyle \left[\frac{1}{30}, 
\frac{7}{30}, 
\frac{11}{30}, 
\frac{13}{30}, 
\frac{17}{30}, 
\frac{19}{30}, 
\frac{23}{30}, 
\frac{29}{30}\right]$} \\ \protect{$\scriptstyle \left[\frac{1}{5}, 
 \frac{1}{4}, 
 \frac{2}{5}, 
 \frac{1}{2}, 
 \frac{3}{5}, 
 \frac{3}{4}, 
 \frac{4}{5}\right]$} \end{array}$ &
            $66$
            \\

            $31$ &
            $\begin{array}{Sc}\protect{$\scriptstyle [30, 1]$}\\ \protect{$\scriptstyle [15, 10, 6]$} \end{array}$ &
            $8$ &
            $\begin{array}{Sc}\protect{$\scriptstyle \left[\frac{1}{30}, 
\frac{7}{30}, 
\frac{11}{30}, 
\frac{13}{30}, 
\frac{17}{30}, 
\frac{19}{30}, 
\frac{23}{30}, 
\frac{29}{30}\right]$} \\ \protect{$\scriptstyle \left[\frac{1}{5}, 
 \frac{1}{3}, 
 \frac{2}{5}, 
 \frac{1}{2}, 
 \frac{3}{5}, 
 \frac{2}{3}, 
 \frac{4}{5}\right]$} \end{array}$ &
            $67$
            \\

            $32$ &
            $\begin{array}{Sc}\protect{$\scriptstyle [15, 2]$}\\ \protect{$\scriptstyle [10, 6, 1]$} \end{array}$ &
            $8$ &
            $\begin{array}{Sc}\protect{$\scriptstyle \left[\frac{1}{15}, 
\frac{2}{15}, 
\frac{4}{15}, 
\frac{7}{15}, 
\frac{8}{15}, 
\frac{11}{15}, 
\frac{13}{15}, 
\frac{14}{15}\right]$} \\ \protect{$\scriptstyle \left[\frac{1}{10}, 
\frac{1}{6}, 
\frac{3}{10}, 
\frac{1}{2}, 
\frac{7}{10}, 
\frac{5}{6}, 
\frac{9}{10}\right]$} \end{array}$ &
            $67$
            \\

            $33$ &
            $\begin{array}{Sc}\protect{$\scriptstyle [15, 7]$}\\ \protect{$\scriptstyle [14, 5, 3]$} \end{array}$ &
            $8$ &
            $\begin{array}{Sc}\protect{$\scriptstyle \left[\frac{1}{15}, 
\frac{2}{15}, 
\frac{4}{15}, 
\frac{7}{15}, 
\frac{8}{15}, 
\frac{11}{15}, 
\frac{13}{15}, 
\frac{14}{15}\right]$} \\ \protect{$\scriptstyle \left[\frac{1}{14}, 
\frac{3}{14}, 
\frac{5}{14}, 
\frac{1}{2}, 
\frac{9}{14}, 
\frac{11}{14}, 
\frac{13}{14}\right]$} \end{array}$ &
            $68$
            \\

            $34$ &
            $\begin{array}{Sc}\protect{$\scriptstyle [30, 5, 3]$}\\ \protect{$\scriptstyle [15, 10, 7, 6]$} \end{array}$ &
            $8$ &
            $\begin{array}{Sc}\protect{$\scriptstyle \left[\frac{1}{30}, 
\frac{7}{30}, 
\frac{11}{30}, 
\frac{13}{30}, 
\frac{17}{30}, 
\frac{19}{30}, 
\frac{23}{30}, 
\frac{29}{30}\right]$} \\ \protect{$\scriptstyle \left[\frac{1}{7}, 
 \frac{2}{7}, 
 \frac{3}{7}, 
 \frac{1}{2}, 
 \frac{4}{7}, 
 \frac{5}{7}, 
 \frac{6}{7}\right]$} \end{array}$ &
            $68$
            \\

            $35$ &
            $\begin{array}{Sc}\protect{$\scriptstyle [30, 5, 3]$}\\ \protect{$\scriptstyle [15, 12, 10, 1]$} \end{array}$ &
            $8$ &
            $\begin{array}{Sc}\protect{$\scriptstyle \left[\frac{1}{30}, 
\frac{7}{30}, 
\frac{11}{30}, 
\frac{13}{30}, 
\frac{17}{30}, 
\frac{19}{30}, 
\frac{23}{30}, 
\frac{29}{30}\right]$} \\ \protect{$\scriptstyle \left[\frac{1}{12}, 
\frac{1}{4}, 
\frac{5}{12}, 
\frac{1}{2}, 
\frac{7}{12}, 
\frac{3}{4}, 
\frac{11}{12}\right]$} \end{array}$ &
            $69$
            \\

            $36$ &
            $\begin{array}{Sc}\protect{$\scriptstyle [15, 6, 1]$}\\ \protect{$\scriptstyle [12, 5, 3, 2]$} \end{array}$ &
            $8$ &
            $\begin{array}{Sc}\protect{$\scriptstyle \left[\frac{1}{15}, 
\frac{2}{15}, 
\frac{4}{15}, 
\frac{7}{15}, 
\frac{8}{15}, 
\frac{11}{15}, 
\frac{13}{15}, 
\frac{14}{15}\right]$} \\ \protect{$\scriptstyle \left[\frac{1}{12}, 
\frac{1}{4}, 
\frac{5}{12}, 
\frac{1}{2}, 
\frac{7}{12}, 
\frac{3}{4}, 
\frac{11}{12}\right]$} \end{array}$ &
            $69$
            \\

            $37$ &
            $\begin{array}{Sc}\protect{$\scriptstyle [15, 1]$}\\ \protect{$\scriptstyle [8, 5, 3]$} \end{array}$ &
            $8$ &
            $\begin{array}{Sc}\protect{$\scriptstyle \left[\frac{1}{15}, 
\frac{2}{15}, 
\frac{4}{15}, 
\frac{7}{15}, 
\frac{8}{15}, 
\frac{11}{15}, 
\frac{13}{15}, 
\frac{14}{15}\right]$} \\ \protect{$\scriptstyle \left[\frac{1}{8}, 
 \frac{1}{4}, 
 \frac{3}{8}, 
 \frac{1}{2}, 
 \frac{5}{8}, 
 \frac{3}{4}, 
 \frac{7}{8}\right]$} \end{array}$ &
            $70$
            \\

            $38$ &
            $\begin{array}{Sc}\protect{$\scriptstyle [30, 5, 3, 2]$}\\ \protect{$\scriptstyle [15, 10, 8, 6, 1]$} \end{array}$ &
            $8$ &
            $\begin{array}{Sc}\protect{$\scriptstyle \left[\frac{1}{30}, 
\frac{7}{30}, 
\frac{11}{30}, 
\frac{13}{30}, 
\frac{17}{30}, 
\frac{19}{30}, 
\frac{23}{30}, 
\frac{29}{30}\right]$} \\ \protect{$\scriptstyle \left[\frac{1}{8}, 
 \frac{1}{4}, 
 \frac{3}{8}, 
 \frac{1}{2}, 
 \frac{5}{8}, 
 \frac{3}{4}, 
 \frac{7}{8}\right]$} \end{array}$ &
            $70$
            \\

            $39$ &
            $\begin{array}{Sc}\protect{$\scriptstyle [20, 3]$}\\ \protect{$\scriptstyle [12, 10, 1]$} \end{array}$ &
            $8$ &
            $\begin{array}{Sc}\protect{$\scriptstyle \left[\frac{1}{20}, 
\frac{3}{20}, 
\frac{7}{20}, 
\frac{9}{20}, 
\frac{11}{20}, 
\frac{13}{20}, 
\frac{17}{20}, 
\frac{19}{20}\right]$} \\ \protect{$\scriptstyle \left[\frac{1}{12}, 
\frac{1}{6}, 
\frac{5}{12}, 
\frac{1}{2}, 
\frac{7}{12}, 
\frac{5}{6}, 
\frac{11}{12}\right]$} \end{array}$ &
            $71$
            \\

            $40$ &
            $\begin{array}{Sc}\protect{$\scriptstyle [20, 6, 1]$}\\ \protect{$\scriptstyle [12, 10, 3, 2]$} \end{array}$ &
            $8$ &
            $\begin{array}{Sc}\protect{$\scriptstyle \left[\frac{1}{20}, 
\frac{3}{20}, 
\frac{7}{20}, 
\frac{9}{20}, 
\frac{11}{20}, 
\frac{13}{20}, 
\frac{17}{20}, 
\frac{19}{20}\right]$} \\ \protect{$\scriptstyle \left[\frac{1}{12}, 
\frac{1}{3}, 
\frac{5}{12}, 
\frac{1}{2}, 
\frac{7}{12}, 
\frac{2}{3}, 
\frac{11}{12}\right]$} \end{array}$ &
            $71$
            \\

            $41$ &
            $\begin{array}{Sc}\protect{$\scriptstyle [20, 1]$}\\ \protect{$\scriptstyle [10, 8, 3]$} \end{array}$ &
            $8$ &
            $\begin{array}{Sc}\protect{$\scriptstyle \left[\frac{1}{20}, 
\frac{3}{20}, 
\frac{7}{20}, 
\frac{9}{20}, 
\frac{11}{20}, 
\frac{13}{20}, 
\frac{17}{20}, 
\frac{19}{20}\right]$} \\ \protect{$\scriptstyle \left[\frac{1}{8}, 
 \frac{1}{3}, 
 \frac{3}{8}, 
 \frac{1}{2}, 
 \frac{5}{8}, 
 \frac{2}{3}, 
 \frac{7}{8}\right]$} \end{array}$ &
            $72$
            \\

            $42$ &
            $\begin{array}{Sc}\protect{$\scriptstyle [20, 3, 2]$}\\ \protect{$\scriptstyle [10, 8, 6, 1]$} \end{array}$ &
            $8$ &
            $\begin{array}{Sc}\protect{$\scriptstyle \left[\frac{1}{20}, 
\frac{3}{20}, 
\frac{7}{20}, 
\frac{9}{20}, 
\frac{11}{20}, 
\frac{13}{20}, 
\frac{17}{20}, 
\frac{19}{20}\right]$} \\ \protect{$\scriptstyle \left[\frac{1}{8}, 
 \frac{1}{6}, 
 \frac{3}{8}, 
 \frac{1}{2}, 
 \frac{5}{8}, 
 \frac{5}{6}, 
 \frac{7}{8}\right]$} \end{array}$ &
            $72$
            \\

            $43$ &
            $\begin{array}{Sc}\protect{$\scriptstyle [20, 1]$}\\ \protect{$\scriptstyle [10, 7, 4]$} \end{array}$ &
            $8$ &
            $\begin{array}{Sc}\protect{$\scriptstyle \left[\frac{1}{20}, 
\frac{3}{20}, 
\frac{7}{20}, 
\frac{9}{20}, 
\frac{11}{20}, 
\frac{13}{20}, 
\frac{17}{20}, 
\frac{19}{20}\right]$} \\ \protect{$\scriptstyle \left[\frac{1}{7}, 
 \frac{2}{7}, 
 \frac{3}{7}, 
 \frac{1}{2}, 
 \frac{4}{7}, 
 \frac{5}{7}, 
 \frac{6}{7}\right]$} \end{array}$ &
            $73$
            \\

            $44$ &
            $\begin{array}{Sc}\protect{$\scriptstyle [20, 7, 2]$}\\ \protect{$\scriptstyle [14, 10, 4, 1]$} \end{array}$ &
            $8$ &
            $\begin{array}{Sc}\protect{$\scriptstyle \left[\frac{1}{20}, 
\frac{3}{20}, 
\frac{7}{20}, 
\frac{9}{20}, 
\frac{11}{20}, 
\frac{13}{20}, 
\frac{17}{20}, 
\frac{19}{20}\right]$} \\ \protect{$\scriptstyle \left[\frac{1}{14}, 
\frac{3}{14}, 
\frac{5}{14}, 
\frac{1}{2}, 
\frac{9}{14}, 
\frac{11}{14}, 
\frac{13}{14}\right]$} \end{array}$ &
            $73$
            \\

            $45$ &
            $\begin{array}{Sc}\protect{$\scriptstyle [20, 3]$}\\ \protect{$\scriptstyle [10, 9, 4]$} \end{array}$ &
            $8$ &
            $\begin{array}{Sc}\protect{$\scriptstyle \left[\frac{1}{20}, 
\frac{3}{20}, 
\frac{7}{20}, 
\frac{9}{20}, 
\frac{11}{20}, 
\frac{13}{20}, 
\frac{17}{20}, 
\frac{19}{20}\right]$} \\ \protect{$\scriptstyle \left[\frac{1}{9}, 
 \frac{2}{9}, 
 \frac{4}{9}, 
 \frac{1}{2}, 
 \frac{5}{9}, 
 \frac{7}{9}, 
 \frac{8}{9}\right]$} \end{array}$ &
            $74$
            \\

            $46$ &
            $\begin{array}{Sc}\protect{$\scriptstyle [20, 9, 6]$}\\ \protect{$\scriptstyle [18, 10, 4, 3]$} \end{array}$ &
            $8$ &
            $\begin{array}{Sc}\protect{$\scriptstyle \left[\frac{1}{20}, 
\frac{3}{20}, 
\frac{7}{20}, 
\frac{9}{20}, 
\frac{11}{20}, 
\frac{13}{20}, 
\frac{17}{20}, 
\frac{19}{20}\right]$} \\ \protect{$\scriptstyle \left[\frac{1}{18}, 
\frac{5}{18}, 
\frac{7}{18}, 
\frac{1}{2}, 
\frac{11}{18}, 
\frac{13}{18}, 
\frac{17}{18}\right]$} \end{array}$ &
            $74$
            \\

            $47$ &
            $\begin{array}{Sc}\protect{$\scriptstyle [24, 1]$}\\ \protect{$\scriptstyle [12, 8, 5]$} \end{array}$ &
            $8$ &
            $\begin{array}{Sc}\protect{$\scriptstyle \left[\frac{1}{24}, 
\frac{5}{24}, 
\frac{7}{24}, 
\frac{11}{24}, 
\frac{13}{24}, 
\frac{17}{24}, 
\frac{19}{24}, 
\frac{23}{24}\right]$} \\ \protect{$\scriptstyle \left[\frac{1}{5}, 
 \frac{1}{4}, 
 \frac{2}{5}, 
 \frac{1}{2}, 
 \frac{3}{5}, 
 \frac{3}{4}, 
 \frac{4}{5}\right]$} \end{array}$ &
            $75$
            \\

            $48$ &
            $\begin{array}{Sc}\protect{$\scriptstyle [24, 5, 2]$}\\ \protect{$\scriptstyle [12, 10, 8, 1]$} \end{array}$ &
            $8$ &
            $\begin{array}{Sc}\protect{$\scriptstyle \left[\frac{1}{24}, 
\frac{5}{24}, 
\frac{7}{24}, 
\frac{11}{24}, 
\frac{13}{24}, 
\frac{17}{24}, 
\frac{19}{24}, 
\frac{23}{24}\right]$} \\ \protect{$\scriptstyle \left[\frac{1}{10}, 
\frac{1}{4}, 
\frac{3}{10}, 
\frac{1}{2}, 
\frac{7}{10}, 
\frac{3}{4}, 
\frac{9}{10}\right]$} \end{array}$ &
            $75$
            \\

            $49$ &
            $\begin{array}{Sc}\protect{$\scriptstyle [24, 4, 1]$}\\ \protect{$\scriptstyle [12, 8, 7, 2]$} \end{array}$ &
            $8$ &
            $\begin{array}{Sc}\protect{$\scriptstyle \left[\frac{1}{24}, 
\frac{5}{24}, 
\frac{7}{24}, 
\frac{11}{24}, 
\frac{13}{24}, 
\frac{17}{24}, 
\frac{19}{24}, 
\frac{23}{24}\right]$} \\ \protect{$\scriptstyle \left[\frac{1}{7}, 
 \frac{2}{7}, 
 \frac{3}{7}, 
 \frac{1}{2}, 
 \frac{4}{7}, 
 \frac{5}{7}, 
 \frac{6}{7}\right]$} \end{array}$ &
            $76$
            \\

            $50$ &
            $\begin{array}{Sc}\protect{$\scriptstyle [24, 7, 4]$}\\ \protect{$\scriptstyle [14, 12, 8, 1]$} \end{array}$ &
            $8$ &
            $\begin{array}{Sc}\protect{$\scriptstyle \left[\frac{1}{24}, 
\frac{5}{24}, 
\frac{7}{24}, 
\frac{11}{24}, 
\frac{13}{24}, 
\frac{17}{24}, 
\frac{19}{24}, 
\frac{23}{24}\right]$} \\ \protect{$\scriptstyle \left[\frac{1}{14}, 
\frac{3}{14}, 
\frac{5}{14}, 
\frac{1}{2}, 
\frac{9}{14}, 
\frac{11}{14}, 
\frac{13}{14}\right]$} \end{array}$ &
            $76$
            \\

            $51$ &
            $\begin{array}{Sc}\protect{$\scriptstyle [24, 4, 3]$}\\ \protect{$\scriptstyle [12, 9, 8, 2]$} \end{array}$ &
            $8$ &
            $\begin{array}{Sc}\protect{$\scriptstyle \left[\frac{1}{24}, 
\frac{5}{24}, 
\frac{7}{24}, 
\frac{11}{24}, 
\frac{13}{24}, 
\frac{17}{24}, 
\frac{19}{24}, 
\frac{23}{24}\right]$} \\ \protect{$\scriptstyle \left[\frac{1}{9}, 
 \frac{2}{9}, 
 \frac{4}{9}, 
 \frac{1}{2}, 
 \frac{5}{9}, 
 \frac{7}{9}, 
 \frac{8}{9}\right]$} \end{array}$ &
            $77$
            \\

            $52$ &
            $\begin{array}{Sc}\protect{$\scriptstyle [24, 9, 6, 4]$}\\ \protect{$\scriptstyle [18, 12, 8, 3, 2]$} \end{array}$ &
            $8$ &
            $\begin{array}{Sc}\protect{$\scriptstyle \left[\frac{1}{24}, 
\frac{5}{24}, 
\frac{7}{24}, 
\frac{11}{24}, 
\frac{13}{24}, 
\frac{17}{24}, 
\frac{19}{24}, 
\frac{23}{24}\right]$} \\ \protect{$\scriptstyle \left[\frac{1}{18}, 
\frac{5}{18}, 
\frac{7}{18}, 
\frac{1}{2}, 
\frac{11}{18}, 
\frac{13}{18}, 
\frac{17}{18}\right]$} \end{array}$ &
            $77$
            \\

\end{longtable}
\bibliographystyle{abbrv}
\bibliography{bib.bib}

\end{document}